\documentclass[leqno,12pt]{article}%
\usepackage{amsmath,amssymb}
\usepackage{dsfont}
\usepackage{amsfonts}%

\newtheorem{Theorem}{\hspace{\parindent}\bf Theorem}[section]
\newtheorem{Lemma}{\hspace{\parindent}\bf Lemma}[section]
\newtheorem{Proposition}{\hspace{\parindent}\bf Proposition}[section]

\renewcommand{\theequation}{\arabic{section}.\arabic{equation}}

\newcommand{\qed}{\hfill$\square$\vspace{0.3cm}}
\begin{document}

\title{\textbf{Nonlocal filtration equations with rough kernels}}
\author{by\\Arturo de Pablo, \\Fernando Quir\'{o}s, and Ana Rodr\'{\i}guez}
\date{}

\maketitle

\begin{center}
\emph{Dedicated to Juan Luis V\'azquez, who has generously shared with us his deep insight on the subject of nonlinear diffusion, on the occasion of his 70th birthday}
\end{center}

\

\begin{abstract}
We study the   nonlinear and nonlocal Cauchy problem
\[
\partial_{t}u+\mathcal{L}\varphi(u)=0 \quad\text{in }\mathbb{R}^{N}\times\mathbb{R}_+,\qquad u(\cdot,0)=u_0,
\]
where $\mathcal{L}$ is a
L\'evy-type nonlocal operator with a kernel having a singularity at the origin as that of the fractional Laplacian. The nonlinearity $\varphi$ is nondecreasing  and continuous, and the initial datum $u_0$ is assumed to be in $L^1(\mathbb{R}^N)$. We prove existence
and uniqueness of weak solutions. For a wide class of nonlinearities, including the porous media case, $\varphi(u)=|u|^{m-1}u$, $m>1$, these solutions turn out to be bounded and   H\"older continuous for $t>0$.   We also describe the large time
behaviour when the nonlinearity resembles a power for $u\approx 0$ and the kernel associated to $\mathcal{L}$ is close at infinity to that of the fractional Laplacian.
\end{abstract}


\vskip 1cm

\noindent{\makebox[1in]\hrulefill}\newline2010 \textit{Mathematics Subject
Classification.} 35R11, 
35S10, 
35B65, 
35K55, 
35B40. 

\newpage

\section{Introduction and main results}

\label{sect-introduction} \setcounter{equation}{0}

We study the   nonlinear and nonlocal Cauchy problem
\begin{equation}
\label{eq:main}
\tag{P}
\left\{
\begin{array}
[c]{ll}%
\partial_{t}u+\mathcal{L}\varphi(u)=0, \qquad & (x,t)\in Q:=\mathbb{R}^{N}\times\mathbb{R}_+,\\[4mm]%
u(x,0) = u_{0}(x), \qquad& x\in\mathbb{R}^{N},
\end{array}
\right.
\end{equation}
with initial data $u_{0}\in L^{1}(\mathbb{R}^{N})$. Sign changes are allowed.   The nonlocal operator
$\mathcal{L}$ is defined formally by
\begin{equation}\label{operatorL}
\mathcal{L} f(x)=\text{P.V.}\int_{\mathbb{R}^{N}}\left( f(x)-f(y)\right) J(x,y)\,dy,
\end{equation}
with a measurable kernel $J$ which is assumed to satisfy
\begin{equation}
\label{eq:kernel}
\tag{$\text{\rm H}_{J}$}
\left\{
\begin{array}{l}
J(x,y)\ge 0, \quad J(x,y)=J(y,x),\\[10pt]
\displaystyle
\frac{\mathds{1}_{\{|x-y|\le
3\}}}{\Lambda|x-y|^{N+\sigma}}\le J(x,y)\le\frac{\Lambda}{|x-y|^{N+\sigma}},
\end{array}
\right.
\qquad  x,\,y\in
\mathbb{R}^{N},\ x\neq y,
\end{equation}
for some constants $\sigma\in(0,2)$ and $\Lambda>0$.   When $J(x,y)=|x-y|^{-(N+\sigma)}$, $\mathcal{L}$ is a multiple of the fractional Laplacian $(-\Delta)^{\sigma/2}$, whose action on smooth functions is well defined and has a pointwise meaning. However, the pointwise expression~\eqref{operatorL} may not have sense for more general kernels in the class that we are considering here, even if $f$ is very smooth.
Hence we have to deal with weak solutions to give sense both to the time derivative and to the nonlocal operator. The  precise definition of a weak solution, in terms of a bilinear form associated to the kernel $J$, is given in Section~\ref{sect-preliminaries}, which is devoted to some preliminaries.

The upper bound in~\eqref{eq:kernel} implies in particular that the operator is of
L\'evy type, $\int_{\mathbb{R}^N}\min(1,| x-  y|^2)J(x,y)\,dy<\infty$ for almost every $x\in\mathbb{R}^N$. Moreover, the singularity on the diagonal $x=y$ is that of the fractional Laplacian.
Thus, $\mathcal{L}$ can be seen as an integro-differential operator of order $\sigma$
with bounded measurable coefficients.
  The bounds in~\eqref{eq:kernel} allow the kernels $J$ to be very oscillating and irregular. That is why they are referred to as \emph{rough} kernels.
Observe that  rapidly decreasing or even compactly supported kernels are permitted. Once they are of L\'evy type, what matters in what follows is their singularity at the origin.

  The linear operator $\partial_t+\mathcal{L}$ is often described in the literature as a nonlocal diffusion operator, since on the one hand it is clear from the integral representation of $\mathcal{L}$ that it is nonlocal, and on the other hand it can be regarded as a diffusion operator, in the sense that solutions to $\partial_t u+\mathcal{L}u=0$ try  to avoid high concentrations.  The same is true for our nonlinear operator $\partial_t+\mathcal{L}\varphi(\cdot)$.

The nonlinearity $\varphi$ is continuous and nondecreasing, and may be assumed without loss of generality to satisfy $\varphi(0)=0$.   The local analogue $\partial_t u-\Delta\varphi(u)=0$ is known as the filtration equation. That is the reason why, by analogy, we label our equation as a \emph{nonlocal filtration equation}.   The typical example is that of powers, $\varphi(s)=|s|^{m-1}s$, which includes both the case of nonlocal porous media, $m>1$, and nonlocal fast diffusion, $0<m<1$.
But we consider also more general functions.

To begin with, in Section~\ref{sect-existence-uniqueness} we prove
existence, uniqueness, and a couple of important properties for bounded weak solutions.
\begin{Theorem}
\label{th:existence} Let $J$ satisfy~\eqref{eq:kernel}, $\varphi$ be continuous and nondecreasing, and  $u_{0}\in L^{1}({\mathbb{R}^{N}})\cap L^{\infty
}({\mathbb{R}^{N}})$.
\begin{itemize}
\item[\rm (a)]  There exists a unique bounded weak solution to the Cauchy
problem~\eqref{eq:main}. It satisfies $\|u(\cdot,t)\|_\infty\le \|u_0\|_\infty$ for every $t>0$.

\item[\rm (b)] If $u$ and $v$ are solutions to problem~\eqref{eq:main}, they satisfy the $T$-contraction property
$$
\int_{\mathbb{R}^N}(u(\cdot,t)-v(\cdot,t))_+\le \int_{\mathbb{R}^N}(u(\cdot,0)-v(\cdot,0))_+\quad\text{for all }t\ge0.
$$

\item[\rm (c)] If moreover $|\varphi(u)|\le C|u|^m$, $m>\frac{(N-\sigma)_+}N$,   and $J(x,y)=\widetilde J(x-y)$,   then
$\int_{\mathbb{R}^N}u(\cdot, t)=\int_{\mathbb{R}^N}u_0$ for all $t\ge0$.
\end{itemize}
\end{Theorem}

\noindent \emph{Remarks. } (i) Existence and uniqueness of bounded \emph{distributional} solutions to~\eqref{eq:main} have been recently obtained in~\cite{Endal-Jakobsen-delTeso} for more  general operators $\mathcal{L}$ than the ones considered here. In the present paper we sacrifice such generality in order to obtain stronger results.

\noindent (ii) We conjecture that conservation of mass is also true in the limit case $m=\frac{N-\sigma}N$ if $N>\sigma$, as it holds when $\mathcal{L}$ is the fractional Laplacian; see~\cite{PQRV2}.

\noindent   (iii) Conservation of mass holds for  kernels satisfying~\eqref{eq:kernel} more general than those considered in paragraph (c), as long as expression~\eqref{operatorL} is well defined for smooth functions; see the beginning of Section~\ref{sect-preliminaries} for some conditions, either on $\sigma$ or on $J$, guaranteeing this fact.

\medskip

We next prove the continuity of bounded solutions when  the nonlinearity satisfies
\begin{equation}
\label{eq:cond.phi}
\tag{$\text{\rm H}_\varphi$}
\varphi\in C^1(\mathbb{R}), \qquad \varphi(0)=0,\qquad \varphi'(s)>0\quad \text{for }s\ne0.
\end{equation}
Notice that $\varphi$ can be degenerate at the level 0. However, we are leaving out nonlinearities which are too degenerate, like the Stefan one, $\varphi(s)=(s-1)_+$,  or singular, like the one corresponding to fast diffusion. If, moreover,
\begin{equation}
\label{eq:cond.phi01}
\tag{$\text{\rm H}_\varphi'$}
\begin{array}{l}
\displaystyle C_1\frac{\varphi(r)}{r}\le \varphi'(s) \le C_2\frac{\varphi(r)}{r}\quad\text{for }0<|r|<1,\ |s|\in(|r|/4,3|r|),
\\[10pt]
\displaystyle
\sup\limits_{[A,B]}\varphi'\le D_M\,\frac{\varphi(B)-\varphi(A)}{B-A}, \qquad\text{if } -M\le A<B\le M,
\end{array}
\end{equation}
for some constants $0<C_1<C_2$ and $D_M>0$, we get H\"older regularity. These conditions control the oscillation of the nonlinearity close to the origin, and are only needed to deal with points at which the equation is degenerate.  They are satisfied for example if
$$
c_1|r|^{m-1}\le\varphi'(r)\le c_2 |r|^{m-1}\quad\text{for }0\le|r|\le1,
$$
for some constants $0<c_1\le c_2,\, m\ge1$.

\begin{Theorem}
\label{th:regularity} Let $J$ and $\varphi$ satisfy respectively~\eqref{eq:kernel} and~\eqref{eq:cond.phi}, and let $u$ be a  bounded weak solution to the
Cauchy problem~\eqref{eq:main}. Then  $u$ is continuous in $Q$. If $\varphi$ satisfies in addition~\eqref{eq:cond.phi01},
then, for all $\tau>0$ there is some $\alpha\in(0,1)$ such that $u\in C^{\alpha}(\mathbb{R}^{N}\times(\tau,\infty))$.
\end{Theorem}

In the proof, performed in Section~\ref{sect-regularity}, we will use De Giorgi's method; see~\cite{deGiorgi}. Thus, we will prove that the oscillation
of the solution in space-time $\sigma$-cylinders of radius $R$,
$$
\mathcal{C}_R=\{|x-x_0|<R,\,|t-t_0|<R^\sigma\}\subset Q,
$$
is
reduced in a fraction of the cylinder $
\mathcal{C}_{\gamma R}$, $\gamma<1$, at least by a constant factor $\varpi_*$. This implies $\sigma$--H\"older continuity,
$$
|u(x,t)-u(x_0,t_0)|\le C\left(|x-x_0|^\alpha+|t-t_0|^{\alpha/\sigma}\right),
$$
with an exponent $\alpha=\log\varpi_*/\log \gamma$.

The control of the oscillation in our nonlocal setting follows the procedure developed in \cite{Caffarelli-Chan-Vasseur} for a linear problem  with a rough kernel (which has applications to certain nonlinear problems), combined with some ideas to deal with the nonlinearity borrowed from~\cite{Athanasopoulos-Caffarelli}.   The operator $\mathcal{L}$   in this latter paper is the fractional Laplacian. This allows to use Caffarelli-Silvestre's extension~\cite{Caffarelli-Silvestre} to transform the problem into a local one. This can not be done for the general kernels that we are considering here. In fact, in the particular case of the fractional Laplacian additional regularity has been obtained~\cite{PQRV3,VPQR}.

  In the linear nonlocal setting, besides~\cite{Caffarelli-Chan-Vasseur}, which uses De Giorgi's technique, we point out the paper~\cite{Felsinger-Kassmann}, where, by means of a different approach based on Moser's work \cite{Moser,Moser-1971}, the authors obtain H\"older regularity with constants that do not depend on the order of differentiability $\sigma\in(0,2)$. See also \cite{Kassmann-2009} for the corresponding elliptic case. It would be interesting to see whether  their method can be adapted to our problem to get rid of the $\sigma$ dependence of the constants.     One can find in the literature other papers dealing with heat kernel estimates and regularity issues for linear parabolic nonlocal problems with rough kernels,  in the framework of Markov jump processes;   see for example \cite{Barlow-Bass-Chen-Kassmann,Bass-Kassmann-Kumagai,Chen,Komatsu} and the references therein.

As for nonlinear   nonlocal problems with rough kernels,   let us mention~\cite{ChangLara-Davila,Serra}, where fully nonlinear nondegenerate  parabolic integro-differential equations are considered.
  Concerning regularity for nonlinear nonlocal equations of porous medium type, besides~\cite{Athanasopoulos-Caffarelli} we have~\cite{Caffarelli-Vazquez-ARMA-2011}, where,  using an approach based on~\cite{Caffarelli-Chan-Vasseur},  the authors prove H\"older regularity for solutions of the so called \emph{porous medium equation with fractional potential pressure},
\begin{equation}\label{eq:otraPME}
\partial_t u=\nabla\cdot(u\nabla p),\qquad p=(-\Delta)^{-\gamma/2}u,\quad \gamma\in(0,2).
\end{equation}

For regularity results for the local filtration analogue to problem~\eqref{eq:main}  we refer to~\cite{dePablo-Vazquez}; see~\cite{Caffarelli-Friedman} for the case of powers.

After this paper was completed, we learned that, at the very same time and independently of us, Bonforte, Figalli and Ros-Oton  proved in \cite{Bonforte-Figalli-RosOton} the H\"older regularity of nonnegative solutions to
the Cauchy-Dirichlet  problem for the fractional porous medium equation $\partial_t u+(-\Delta)^{\sigma/2} u^m$, $m>1$, in a bounded domain. The authors also indicate how the result could be extended to general unbounded domains in $\mathbb{R}^N$ for equations of the form~\eqref{eq:main}. The method of proof in that paper is completely different from ours and does not apply to solutions with sign changes.

The assumption \lq\lq$u$ bounded'' in the previous regularity result is not a big restriction as we see now. Indeed, if the kernel $J$ satisfies the stronger assumption
\begin{equation}\label{eq:kernel2}
\tag{$\text{\rm H}_{J}'$}
\frac{1}{\Lambda|x-y|^{N+\sigma}}\le J(x,y)\le\frac{\Lambda}{|x-y|^{N+\sigma}},\qquad x,\,y\in\mathbb{R}^N,\;x\ne y,
\end{equation}
besides~\eqref{eq:kernel}, the natural energy  associated to the operator $\mathcal{L}$ is in fact equivalent to the fractional Sobolev energy $\|(-\Delta)^{\sigma/4}\|_2^2$; see Section~\ref{sect-preliminaries}. Then it is possible to get an $L^1$--$L^\infty$ smoothing effect repeating the Moser-like arguments in~\cite{VPQR}. This allows in addition to get an existence and uniqueness result for initial values in $L^1(\mathbb{R}^N)$ by approximation.

\begin{Theorem}
\label{th:smoothing}
Let $J$ satisfy \eqref{eq:kernel2} and $\varphi\in C^1(\mathbb{R}\setminus\{0\})$ be such that $\varphi'(s)\ge C |s|^{m-1}$ for
some  $m>\frac{(N-\sigma)_+}{N}$. If $u_0\in L^1({\mathbb{R}^N})$, and $N>\sigma$, there exists a unique weak solution to the Cauchy problem~\eqref{eq:main} which is bounded in $\mathbb{R}^N\times(\tau,\infty)$ for all  $\tau>0$. This solution moreover satisfies
\begin{equation}
\label{eq:L-inf-p3}
\|u(\cdot,t)\|_{\infty}\le C_{1}t^{-\gamma }\|u_{0}
\|_{1}^{\delta},
\end{equation}
with $\gamma=\frac{N}{N(m-1)+{\sigma}}$ and $\delta=\frac{\sigma \gamma}N$,
the constant $C_{1}$ depending on $m,\,\sigma, C,$ and $N$.

If $N=1\le\sigma<2$ the result is still valid if we assume in addition $\varphi'(s)\le \widetilde C |s|^{m-1}$.
\end{Theorem}

We  next turn our attention, in Section~\ref{sect-asymptotic-behaviour}, to the asymptotic behaviour of the solutions when the operator $\mathcal{L}$ behaves in some sense as $(-\Delta)^{\sigma/2}$ and the constitutive function $\varphi$ behaves as a power in a neighbourhood of the origin. To be more precise,  we assume that
\begin{align}
\label{eq:J.symmetric}
&J(x,y)= \widetilde J(z),\quad z=|x-y|;
\\
\label{eq:cond.J}
&\lim_{|z|\to\infty}|z|^{N+\sigma}\widetilde J(z)=\mu>0;
\\
\label{eq:cond.phi.prima}
&\displaystyle \lim_{u\to0}|u|^{1-m}\varphi'(u)=a>0\text{ for some }m\ge1.
\end{align}
Under these conditions, we will prove that
the solution behaves for large times as the solution
$B=B_M$ to
\begin{equation}
\label{eq:FPME}
\left\{
\begin{array}
[c]{ll}%
\partial_{t}B+\frac{a\mu}{m\mu_{N,\sigma}}(-\Delta)^{\sigma/2}\left(|B|^{m-1}B\right)=0\quad&\text{in }Q,\\[4mm]%
B(\cdot,0) = M\delta\quad&\text{in }\mathbb{R}^{N},
\end{array}
\right.
\end{equation}
where $\delta$ is the unit Dirac mass placed at the origin.
The constant $\mu_{N,\sigma}$, which is explicit, appears  as a normalization constant in the definition of $(-\Delta)^{\sigma/2}$. Since the mass is preserved in the evolution, see~Theorem~\ref{th:existence}, then necessarily $M=\int_{\mathbb{R}^N}u_0$.
The equation in~\eqref{eq:FPME} has been analysed
in~\cite{PQRV,PQRV2} when the initial data are integrable, and in \cite{Vazquez} when they are  non-negative Radon measures.
The function $B_M$, obtained in the latter reference,  is called a fundamental, or  Barenblatt,
solution; see~\cite{Blumenthal-Getoor} for the linear case $m=1$.  It has a definite sign, given by $M$, and a self-similar
structure,
\begin{equation}
\label{eq:Barenblatt}
B_M(x,t)=t^{-\alpha}Z(xt^{-\beta}), \qquad
\alpha=\frac N{N(m-1)+\sigma}, \quad \beta=\frac1{N(m-1)+\sigma}.
\end{equation}

Compactness will follow from the H\"older estimates provided by~Theorem~\ref{th:regularity}, which indeed hold, thanks to~\eqref{eq:cond.phi.prima}.

\begin{Theorem}
\label{th:asymptotic.behaviour}
Let $J$ and $\varphi$ satisfy respectively~\eqref{eq:kernel} and~\eqref{eq:cond.phi} and also \eqref{eq:J.symmetric}--\eqref{eq:cond.phi.prima}. Let $u$ be a  bounded solution to~\eqref{eq:main}, where $\varphi\in C^{1,\gamma}(\mathbb{R})$ for some $\gamma\in(0,1)$, and $M=\int_{\mathbb{R}^N}u_0$. Then
\begin{equation}
\label{eq:asymptotic.behaviour}
t^{\frac N{N(m-1)+\sigma}}\|u(\cdot,t)-B_M(\cdot,t)\|_\infty\to0\quad\text{as }t\to\infty.
\end{equation}
\end{Theorem}

The proof of~\eqref{eq:asymptotic.behaviour} uses the existence of a solution with certain properties for a linear (dual) problem with coefficients in nondivergence form,  which corresponds to a nonsymmetric kernel.  This problem, which has independent interest, will be studied in the Appendix.

The corresponding result for the   (local) case in which $\mathcal{L}=-\Delta$   was first obtained in~\cite{Kamin}; see also~\cite{Biler-Dolbeault-Esteban,Carrillo-diFrancesco-Toscani}. The only known result up to now for the non-local case is given in~\cite{Vazquez}, where the author obtains the asymptotic behaviour for non-negative solutions for the problem with $\varphi(u)=u^m$, $m> \frac{(N-\sigma)_+}{N}$, and $\mathcal{L}=(-\Delta)^{\sigma/2}$.   Let us also mention the work~\cite{Caffarelli-Vazquez-DCDS-2011} on the asymptotic behaviour of solutions to the porous medium equation with fractional potential pressure~\eqref{eq:otraPME}. The large time behaviour is given again by a Barenblatt type solution, which was constructed in~\cite{Biler-Imbert-Karch-2015}, and which turns out to be related to the Barenblatt solution of problem~\eqref{eq:main}, as proved in~\cite{Stan-delTeso-Vazquez-2015}.

\medskip

\noindent\emph{Remarks. }  (i) Notice that the precise behaviour of $ \widetilde J $ is only needed at infinity. This is due to the fact that mass goes to zero in compact sets. The behaviour far from infinity is only needed to obtain compactness, via Theorem~\ref{th:regularity}.

\noindent (ii) The boundedness of the solution is not a restriction for a wide class of  nonlinearities; see Theorem~\ref{th:smoothing}.

\noindent (iii) If $M=0$, the result does  not give a non-trivial asymptotic profile, but only that $u=o\left(t^{-\frac N{N(m-1)+\sigma}}\right)$, which is nevertheless better that what the smoothing effect gives; see~\eqref{eq:L-inf-p3}. A  nontrivial self-similar asymptotic behaviour is still expected; see~\cite{Kamin-Vazquez} for the local case with power nonlinearities, where a limit dipole solution is obtained for $N=1$. This case will be treated elsewhere.


\section{Preliminaries}

\label{sect-preliminaries} \setcounter{equation}{0}

In this section we establish the notion of weak solution to problem~\eqref{eq:main} and fix the required functional
framework.

  As mentioned in the Introduction, expression~\eqref{operatorL} is only formal and may not make sense for the general kernels that we are considering here, even for smooth functions. The validity of~\eqref{operatorL} is guaranteed only for $\sigma< 1$ and functions $f$ in $C^{\sigma+\varepsilon}$ that do not grow too much at infinity.  However, we notice  that if we assume the additional condition
\begin{equation}
\label{eq:cond.J.poinwise.L}
J(x,x+y)=J(x,x-y),
\end{equation}
then the operator has a pointwise expression, in terms of second differences, even for $1\le\sigma<2$, for regular  enough, $C^{1,\sigma-1+\varepsilon}$, functions,
\begin{equation}\label{operatorL2}
\mathcal{L} f(x)=-\frac12\int_{\mathbb{R}^{N}}\left( f(x+y)+f(x-y)-2f(x)\right) J(x,x-y)\,dy.
\end{equation}
In the case $J(x,y)=\widetilde J(x-y)$,  condition~\eqref{eq:cond.J.poinwise.L} follows from the symmetry of the kernel. We remark that, except in Section~\ref{sect-asymptotic-behaviour}, where we assume that $J(x,y)=\widetilde J(|x-y|)$,  we will not impose~\eqref{eq:cond.J.poinwise.L}. In any case, even if~\eqref{eq:cond.J.poinwise.L} holds, solutions to problem~\eqref{eq:main} need not be classical and we have to consider weak solutions, and in particular a weak definition of the operator.

In order to define the action of the operator $\mathcal{L}$ in a weak sense we  consider the bilinear form
(nonlocal interaction energy)%
\begin{equation*}
\label{quadratic-form}\mathcal{E}(f,g)=\frac12\int_{\mathbb{R}^{N}}%
\int_{\mathbb{R}^{N}}(f(x)-f(y))(g(x)-g(y))J(x,y)\,dxdy,
\end{equation*}
and the quadratic form $\overline{\mathcal{E}}(f)=\mathcal{E}(f,f)$.   For kernels satisfying the symmetry condition~\eqref{eq:cond.J.poinwise.L} and functions $f,\, g\in C^2_0(\mathbb{R}^N)$ we have
$$
\langle \mathcal{L}f,g\rangle=\mathcal{E}(f,g);
$$
see~\cite{Schilling-Uemura}.
The bilinear form   $\mathcal{E}$ is well defined for more general kernels, not necessarily satisfying~\eqref{eq:cond.J.poinwise.L},
and for functions in the space  $\dot{\mathcal{H}}_{\mathcal{L}}(\mathbb{R}^{N})$, which is the closure of $C_0^\infty(\mathbb{R}^N)$ with the seminorm associated to the quadratic form $\overline{\mathcal{E}}$. We also define
\[
\label{def-HL}\mathcal{H}_{\mathcal{L}}(\mathbb{R}^{N})=\{f\in L^{2}%
(\mathbb{R}^{N})\,:\,\overline{\mathcal{E}}(f)<\infty\}.
\]
When $J(x,y)=|x-y|^{-N-\sigma}$ for some $0<\sigma<2$, the operator reduces to a multiple of the fractional Laplacian $(-\Delta)^{\sigma/2}$. It is clear then from \eqref{eq:kernel} that the space $\mathcal{H}_{\mathcal{L}}(\mathbb{R}^{N})$ coincides with the
fractional Sobolev space
$$
H^{\sigma/2}(\mathbb{R}^{N})=\{f\in L^{2}%
(\mathbb{R}^{N})\,:\,(-\Delta)^{\sigma/4}\in L^{2}%
(\mathbb{R}^{N})\}.
$$
Actually,
\eqref{eq:kernel} implies
\begin{equation}
\label{estimate-norm-HL}c_{1}\,\overline{\mathcal{E}}(f) \le\|(-\Delta
)^{\sigma/4}f\|_{2}^{2}\le c_{2}\left(\|f\|_{2}^{2}+
\overline{\mathcal{E}}(f)\right),
\end{equation}
where $c_1,\,c_2$ depend only on $N,\sigma,\Lambda$.
If we assume \eqref{eq:kernel2}, we get the stronger result $\overline
{\mathcal{E}}(f)\sim\|(-\Delta)^{\sigma/4}f\|_{2}^{2}$.

We recall also the inclusions $H^{\sigma/2}(\mathbb{R}^{N})\subset L^{\frac{2N}{N-\sigma}}(\mathbb{R}^{N})$  if $N>\sigma$, and
$H^{\sigma/2}(\mathbb{R})\subset L^q(\mathbb{R})$ for every $q\ge2$ if $1\le\sigma<2$. More precisely, we have the
Hardy-Littlewood-Sobolev inequality \cite{Hardy-Littlewood,Sobolev},
\begin{equation}
\label{HLS}\|(-\Delta)^{\sigma/4}f\|_{2}\ge c\|f\|_{\frac{2N}{N-\sigma}},\qquad N>\sigma,
\end{equation}
and the Nash-Gagliardo-Nirenberg inequality \cite{PQRV2},
\begin{equation}
\label{NGN}\|(-\Delta)^{\sigma/4}f\|^2_{2}\|f\|^p_{p}\ge c\|f\|^{p+2}_{\frac{N(p+2)}{2N-\sigma}},\qquad N\ge1,\;0<\sigma<2,\;p\ge1.
\end{equation}
These two inequalities combined with the upper estimate in~\eqref{estimate-norm-HL} yield useful inclusions for functions in $\mathcal{H}_{\mathcal{L}}$.

When dealing with bounded domains $\Omega\subset\mathbb{R}^N$ we consider the operator acting on functions vanishing outside $\Omega$. The corresponding Sobolev type space is $\mathcal{H}_{\mathcal{L},0}(\Omega)$, defined as the completion of $C_0^\infty(\Omega)$ with the norm  given by $\overline{\mathcal{E}}^{1/2}$.
 Functions in this space satisfy a Poincar\'e inequality \cite{Andreu-Mazon-Rossi-Toledo},
\begin{equation}
  \label{poincare}
  \overline{\mathcal{E}}(f)\ge c(\Omega)\|f\|_2^2.
\end{equation}

We now define the concept of \emph{weak solution} to the Cauchy
problem~\eqref{eq:main}, a function $u\in C([0,\infty): L^{1}(\mathbb{R}%
^{N}))$ with $\varphi(u) \in L^2_{\rm
loc}((0,\infty):\dot{\mathcal{H}}_{\mathcal{L}}(\mathbb{R}^{N}))$ such that
\begin{equation}
\label{weak-nonlocal}\displaystyle \int_{0}^{\infty}\int_{\mathbb{R}^{N}%
}u\partial_{t} \zeta-\int_{0}^{\infty}\mathcal{E}(\varphi(u),\zeta
)=0
\end{equation}
for every $\zeta\in C_{c}^{\infty}(Q)$, and taking the initial datum
$u(\cdot,0)=u_{0}$ almost everywhere.

One of the tools needed in the following sections is a generalized
Stroock-Varopoulos inequality; see~\cite{Brandle-dePablo,Varopoulos}.

\begin{Proposition}
\label{pro:SV-J}  If $F,\,G$ are two functions such that $F(u),\,G(u)\in
\dot{\mathcal{H}}_{\mathcal{L}}(\mathbb{R}^{N})$, then
\begin{equation}
\label{eq:SV-J}\overline{\mathcal{E}}(H(u))\le\mathcal{E}(F(u),G(u)),
\end{equation}
if $(H^\prime)^2 \le F^{\prime}
G^{\prime}$.
\end{Proposition}


\section{Existence and uniqueness. Proof of Theorem~\ref{th:existence}}

\label{sect-existence-uniqueness} \setcounter{equation}{0}

In order to perform the existence proof we rewrite the
problem in the equivalent form
\begin{equation}
\label{eq:main-beta}
\tag{$\text{P}_\beta$}
\partial_{t}\beta(w)+\mathcal{L}w=0,
\end{equation}
where $w=\varphi(u)$ and $\beta=\varphi^{-1}$.

We  construct solutions by means of Crandall-Liggett's
Theorem~\cite{Crandall-Liggett}, which is based on an implicit in time
discretization. Hence, we  have to deal with the elliptic problem
\begin{equation}
\label{elliptic}\beta(v)+\mathcal{L}v=g\qquad\mbox{in } \mathbb{R}^{N},%
\end{equation}
with $g\in L^1(\mathbb{R}^{N})\cap L^\infty(\mathbb{R}^{N})$.
To show existence of a weak solution for this problem we approximate the space
$\mathbb{R}^{N}$ by finite balls $B_{R}$. Thus, we look for a weak
solution $v=v_{R}\in \mathcal{H}_{\mathcal{L},0}(B_{R})$ to the problem,
\begin{equation}
\label{elliptic-bdd}
\beta(v)+\mathcal{L}v=g\quad\text{in } B_{R},\qquad
v=0 \quad\text{in } B_{R}^{c},
\end{equation}
that is,
\begin{equation*}
\label{weak.elliptic.bdd}{\mathcal{E}}(v,\zeta)+\int_{B_{R}}\beta(v)\zeta
-\int_{B_{R}}g\zeta=0,
\end{equation*}
for every test function $\zeta\in \mathcal{H}_{\mathcal{L},0}(B_{R})$. Existence is obtained in a standard way by
minimizing the functional
\[
J(v)=\frac12\,\overline{\mathcal{E}}(v)+\int_{B_{R}}\Theta(v)- \int_{B_{R}}vg
\]
in $\mathcal{H}_{\mathcal{L},0}(B_{R})$, where $\Theta^{\prime}=\beta$. This functional is coercive in $\mathcal{H}_{\mathcal{L},0}(B_R)$. Indeed,  using H\"older's inequality, we have, for every $\varepsilon>0$,
$$
\left|\int_{B_R} vg\right|\le
\|v\|_{\frac{2N}{N-\sigma}}\,\|g\|_{\frac{2N}{N+\sigma}}\le
\varepsilon\|v\|^2_{\frac{2N}{N-\sigma}}+\frac1{4\varepsilon}\|g\|_{\frac{2N}{N+\sigma}}^2.
$$
Thus, Hardy-Littlewood-Sobolev inequality \eqref{HLS}, together with Poincar\'e inequality \eqref{poincare}, implies, if $N>\sigma$,
\begin{equation*}
J(v)\ge C_1 \overline{\mathcal{E}}(v)-C_2 .
\label{coercivo}
\end{equation*}
For $N=1\le\sigma<2$ we use the Nash-Gagliardo-Nirenberg
inequality \eqref{NGN} instead.

We have thus obtained a weak solution $v_R$ to~\eqref{elliptic-bdd}.
On the other hand, given two data $g_1$ and $g_2$, the
corresponding weak solutions satisfy the $T$-contraction property
\begin{equation*}\label{eq:-contraction-bdd}
\int_{B_R}\left(\beta(v_{R,1})-\beta(v_{R,2})\right)_+\le
\int_{B_R}\left(g_1- g_2\right)_+.
\end{equation*}
In particular, $\|\beta(v_R)\|_{L^1(B_R)}\le\|g\|_{L^1(B_R)}$ and  $\|\beta(v_R)\|_{L^\infty(B_R)}\le\|g\|_{L^\infty(B_R)}$.
It is then easy to prove that the monotone  limit
$v=\lim_{R\to\infty}v_R$ is a weak solution to
problem \eqref{elliptic}.
The $T$-contractivity property
also holds in the limit. Moreover,
$\|\beta(v)\|_{L^\infty(\mathbb{R}^N)}\le\|g\|_{L^\infty(\mathbb{R}^N)}$
and  $\|\beta(v)\|_{L^1(\mathbb{R}^N)}\le\|g\|_{L^1(\mathbb{R}^N)}$.

Now, using Crandall-Liggett's
Theorem we obtain the existence of a unique mild solution
$w$ to the evolution problem \eqref{eq:main-beta}. It is moreover a weak solution since it lies in the energy
space. This is checked using the same technique as in~\cite{PQRV},
which yields, taking $\Phi'=\varphi$,
\begin{equation*}\label{L2grad}
\int_0^T\overline{\mathcal{E}}(w)dt\le \int_{\mathbb{R}^N}\Phi(u_0)\le \|u_0\|_1 \|\varphi(u_0)\|_\infty\quad \text{for every }T>0.
\end{equation*}

Uniqueness follows by the standard argument due to Oleinik et al.~\cite{Oleinik-Kalashnikov-Czou}; see \cite{VPQR}. The parabolic  $T$-contraction can be deduced from its elliptic counterpart.

In order to complete the proof  we show now that the conservation of mass is true if $|\varphi(u)|\le C|u|^{m}$ above the critical exponent $\frac{(N-\sigma)_+}N$,   and  $J(x,y)=\widetilde J(x-y)$.
We  adapt the technique used in the local case. Take a nonnegative
non-increasing smooth cut-off function $\psi(s)$ such that $\psi(s)=1$ for
$0\le s\le1$, $\psi(s)=0$ for $s\ge2$, and define
$\phi_R(x)=\psi(|x|/R)$.   Since $\mathcal{L}\phi_R$ is well defined under our assumptions on $J$,  we obtain, for every $t>0$,
\begin{equation}
\displaystyle\int_{\mathbb{R}^N}
u(\cdot,t)\phi_R-\int_{\mathbb{R}^N}u_0\phi_R=-\int_0^t\mathcal{E}(\varphi(u),\phi_R)=-\int_0^t\int_{{\mathbb{R}^N}}\varphi(u)\mathcal{L}\phi_R.
\label{eq:mass-cons}
\end{equation}
On the other hand,
the radial cut-off function $\phi_R$ has the scaling property
\begin{equation*}
\label{varphi1}
\mathcal{L}\phi_R(x)=R^{-\sigma}\widetilde{\mathcal{L}}\phi_1(x/R),
\end{equation*}
where $\widetilde{\mathcal{L}}$ is another nonlocal operator satisfying the same properties as $\mathcal{L}$. In particular $\widetilde{\mathcal{L}}\phi_1\in L^1(\mathbb{R}^N)\cap
L^\infty(\mathbb{R}^N)$.  This implies $\|\mathcal{L}\phi_R\|_q\le CR^{-\sigma+N/q}$ for every $1\le q\le\infty$.
Then, if we apply H\"older's inequality with $p=\max\{1,1/m\}$ to
the right-hand side of \eqref{eq:mass-cons}, and use the above
property, together with the estimate $|\varphi(u)|\le C|u|^m$, we get

$$
\begin{array}{rl}
\displaystyle\left|\int_{\mathbb{R}^N}u(\cdot,t)\phi_R-\int_{\mathbb{R}^N}u_0\phi_R\right|&\displaystyle\le t\|u_0\|_{\infty}^{m-1/p}
\|u_0\|_{1}^{1/p}\|{\mathcal{L}}\phi_R\|_{p/(p-1)} \\
&\displaystyle\le t C R^{-\sigma+N(p-1)/p}\|u_0\|_{\infty}^{m-1/p}
\|u_0\|_{1}^{1/p}.
\end{array}
$$
The result follows letting $R$ go to infinity, since the exponent of $R$ is negative precisely for
$m>\frac{N-\sigma}N$.


\section{Regularity. Proof of Theorem~\ref{th:regularity}}

\label{sect-regularity} \setcounter{equation}{0}

As in the previous section, it is convenient to work with equation~\eqref{eq:main-beta}, so that the nonlocal term is linear, the nonlinearity being confined to the time derivative. In the course of the proof we will need to establish some estimates for the solutions of
 \begin{equation}
\label{eq:main-vartheta}
\tag{$\text{P}_\vartheta$}
\partial_{t}\vartheta(w)+\mathcal{K}w=0,
\end{equation}
for different functions $\vartheta$ and operators $\mathcal{K}$ related, respectively, to our original function $\beta=\varphi^{-1}$ and our original nonlocal operator $\mathcal{L}$. To be more precise, $\vartheta$ will have the form $\vartheta(s)=a\beta(bs+c)$ for some $a,\, b>0$. Hence, in the sequel we always assume without further mention that
\begin{equation*}
\label{eq:conditions.vartheta}
\vartheta\in C(\mathbb{R})\cap C^1(\mathbb{R}\setminus\{s_0\}),\quad \vartheta'(s)>0\quad\text{for }s\neq s_0,\qquad \text{for some }s_0\in\mathbb{R}.
\end{equation*}
As for the operator $\mathcal{K}$, its kernel will always satisfy~\eqref{eq:kernel}.

For any given Lipschitz function $\psi$, we define the functional
$$
\mathcal{B}_{\psi}(v)=\int_{0}^{(v-\psi)_{+}}\vartheta^{\prime}(s+\psi)s\,ds
$$
The first step of the regularity argument is to obtain, by using the equation, an energy estimate associated to $\mathcal{B}_\psi$. The quadratic form $\overline{\mathcal{E}}$ and the bilinear
form $\mathcal{E}$ always refer to the operator $\mathcal{K}$ being considered.

\begin{Lemma}
Let $\psi\in C^{0,1}(\mathbb{R}^{N})$ satisfy $\int_{\{|x-y|>1\}}|\psi(x)-\psi(y)|J(x,y)\,dy<C<\infty$ for every $x\in\mathbb{R}^N$, and $w$ be a weak solution to~\eqref{eq:main-vartheta} in some finite time interval $I$ including $(t_1,t_2)$. Then,
\begin{equation}
\label{energy3}
\begin{array}
[c]{l}%
\displaystyle\int_{\mathbb{R}^{N}}\mathcal{B}_{\psi}(w)(x,t_{2})\,dx+ \int_{t_{1}%
}^{t_{2}}\overline{\mathcal{E}}((w-\psi)_{+})(t)\,dt\le\int_{\mathbb{R}^{N}%
}\mathcal{B}_{\psi}(w)(x,t_{1})\,dx\\[10pt]%
\quad\displaystyle +C\int_{t_{1}}^{t_{2}}\left( \int_{\mathbb{R}^{N}}
(w-\psi)_{+}(x,t)+\mathds{1}_{\{w(x,t)> \psi(x)\}}\right) \,dx dt.
\end{array}
\end{equation}
\end{Lemma}

\noindent\emph{Proof. }
If we multiply equation~\eqref{eq:main-vartheta} by the function $\zeta=(w-\psi)_{+}$, we formally get
\begin{equation}
\label{weak-w}\left. \int_{\mathbb{R}^{N}}\mathcal{B}_{\psi}(w(x,t))\,dx\right| _{t_{1}%
}^{t_{2}}+\int_{t_{1}}^{t_{2}}\mathcal{E}(w,(w-\psi)_{+})(t)\,dt=0.
\end{equation}
Though $w$
is not differentiable in time almost everywhere, a regularization procedure in the weak
formulation using some Steklov averages, following an idea from
\cite{Aronson-Serrin}, allows to bypass this difficulty. In fact, it suffices to show that
$\partial_t\mathcal{B}_{\psi}(w)\in L^{2}_{\mathrm{loc}}(I:L^{2}(\mathbb{R}^{N}))$.

 For any $g\in L^{1}(\mathbb{R}^N\times I)$ we
define the Steklov average
\[
g^{h}(x,t)=\frac1h\int_{t}^{t+h}g(x,s)\,ds.
\]
We see that almost everywhere we have
\[
\partial_{t}g^{h}(x,t)=\delta^{h}g(x,t):=\frac{g(x,t+h)-g(x,t)}{h}\,.
\]
Since $\partial_{t}\vartheta(w)^{h}\in L^{1}(\mathbb{R}^{N}\times I)$, we can write the
weak formulation~\eqref{weak-nonlocal} in the form
\[
\int_{0}^{\infty}\int_{\mathbb{R}^{N}}\partial_{t}\vartheta(w)^{h}\zeta
=\int_{0}^{\infty}\mathcal{E}(w^{h},\zeta).
\]
To simplify we perform the calculations with $\psi=0$.
We take $\zeta=(\chi\partial_{t}w^{h})^{-h}$ as test function, where $\chi\in
C_{0}^{\infty}(I)$, $0\le\chi\le1$, $\chi(t)=1$ for $t\in[t_{1}%
,t_{2}]$, is a cut-off function. Using the ``integration by parts'' formulae
$\displaystyle\int_{0}^{\infty}f\delta^{h}g=-\int_{0}^{\infty}g\delta^{-h}f$,
and $\displaystyle \int_{0}^{\infty}\mathcal{E}(f,g^{h})=\int_{0}^{\infty
}\mathcal{E}(f^{-h},g)$, the above identity becomes
\[
\int_{0}^{\infty}\int_{\mathbb{R}^{N}}\chi\partial_{t}\vartheta(w)^{h}%
\,\partial_{t}w^{h} =\frac12\int_{0}^{\infty}\chi\partial_{t}%
\mathcal{E}(w^{h},w^{h}) =-\frac12\int_{0}^{\infty}\partial_{t}%
\chi\mathcal{E}(w^{h},w^{h}).
\]
At this point we observe that the same calculus inequality used in
\cite{Brandle-dePablo} to prove \eqref{eq:SV-J} allows to show that $\delta^{h} \vartheta(w)\,\delta^{h} w\ge(\delta^{h}(\ell(w)))^{2}$,
where $(\ell^\prime)^2=\vartheta^{\prime}$. We therefore get
\[
\int_{t_{1}}^{t_{2}}\int_{\mathbb{R}^{N}}(\delta^{h}(\ell(w)))^{2}\le
c\int_{0}^{\infty}|\partial_{t}\chi^{\prime}|\,\mathcal{E}(w^{h},w^{h})\le
c.
\]
On the other hand, $|\delta^{h}B(w)| \le|\sqrt{\vartheta^{\prime}(w)}w\delta
^{h}\ell(w)|$, so $\delta^{h}B(w)\in L^{2}(I:L^{2}(\mathbb{R}^{N}))$ provided $\sqrt{\vartheta^{\prime}(w)}w\in L^\infty(\mathbb{R}^N\times I)$, and we end by passing to the limit $h\to0$.

Once we have~\eqref{weak-w}, the energy estimate is  obtained proceeding as in~\cite{Caffarelli-Chan-Vasseur}.
\qed

A consequence of this energy estimate is obtained using the properties of $\vartheta$ and $\mathcal{E}$. If  $\ell=\inf_{\{w\ge\psi\}}w\ge0$ and $M=\sup_{\{w\ge\psi\}}w<\infty$, we have
\begin{equation}
\label{estimate-B}\Lambda_{1}(w-\psi)_+^{2}\le
\mathcal{B}_{\psi}(w)\le \Lambda_{2} (w-\psi)_+,
\end{equation}
where
\begin{equation}
\label{eq:lambdas}
\Lambda_1=\frac12\inf_{\ell\le s\le M}\vartheta'(s),\qquad \Lambda_2=\vartheta(M)-\vartheta(\ell).
\end{equation}
Therefore, using \eqref{estimate-norm-HL}, the energy estimate~\eqref{energy3} yields
\begin{equation}\label{eq:newenergy}
\begin{array}
[c]{l}%
\displaystyle\Lambda_{1}\int_{\mathbb{R}^{N}}(w-\psi)_+^{2}(x,t_2)\,dx+ c\int_{t_{1}%
}^{t_2}\int_{\mathbb{R}^{N}}|(-\Delta)^{\sigma/4}((w-\psi)_{+})(x,t)|^2\,dxdt\\[10pt]%
\quad\displaystyle\le\Lambda_{2}\int_{\mathbb{R}^{N}%
}(w-\psi)_+(x,t_{1})\,dx\\[10pt]
\qquad\displaystyle+C\int_{t_{1}}^{t_2} \int_{\mathbb{R}^{N}}
\Big((w-\psi)^2_{+}(x,t)+(w-\psi)_{+}(x,t)+\mathds{1}_{\{w(x,t)> \psi(x)\}}\Big) \,dx dt.
\end{array}
\end{equation}
This is a kind of \lq\lq Anti-Sobolev inequality'', controlling the energy in terms of the size of the solution.

The next step is to obtain a first De Giorgi type oscillation reduction lemma: if $w$ is mostly negative in space-time measure in a certain parabolic cylinder, then the supremum goes down if we restrict to a smaller nested cylinder.
Due to the nonlocal character of the operator, it is necessary to have some control of the far away behaviour of the solution. This is done, as in \cite{Caffarelli-Chan-Vasseur}, via a barrier function. In order to simplify our approach we work with normalized cylinders. The general case is treated by scaling.

\noindent\emph{Notation. } $\Gamma_R=B_R(0)\times[-R^\sigma,0]$.

%

\begin{Lemma}
\label{lem:first.DG.lemma}
If $\vartheta$ satisfies
\begin{equation}
\label{eq:delta.vartheta}
\delta(\vartheta):=\frac{\inf\limits_{0\le s\le 2}\vartheta'(s)}{1+\vartheta(2)-\vartheta(0)}>0,
\end{equation}
there is a constant $c>0$ such that if $w:\mathbb{R}^N\times(-2,0)\to \mathbb{R}$ is a  weak solution to equation~\eqref{eq:main-vartheta} satisfying
\begin{align}
\label{eq:primerapsi}
&w(x,t)\le 1+(|x|^{\sigma/4}-1)_{+}\quad\text{in } \mathbb{R}^N\times(-2,0),\\[10pt]
\label{eq:critical.condition}
&|\{w>0\}\cap\Gamma_2|\le c\delta(\vartheta)^{1+N/\sigma},
\end{align}
then
$$
w(x,t)\le \frac12\quad\text{if }(x,t)\in \Gamma_1.
$$
\end{Lemma}

Once \eqref{estimate-B} is true, we can perform the same proof of \cite[Lemma 3.1]{Caffarelli-Chan-Vasseur}. Nevertheless, on the one hand we have to pursue the constants $\Lambda_i$ in~\eqref{eq:lambdas}, to see how the nonlinearity $\vartheta$ affects the result, in the spirit of~\cite{Athanasopoulos-Caffarelli}. This gives the precise value of $\delta(\vartheta)$ in~\eqref{eq:delta.vartheta}. On the other hand, the proof performed in~\cite{Caffarelli-Chan-Vasseur} works only for $N>\sigma$ since Hardy-Littlewood-Sobolev inequality~\eqref{HLS} is used. We complement the result for $N=1\le\sigma<2$ by using  Nash-Gagliardo-Nirenberg inequality \eqref{NGN}.

\noindent\emph{Proof. }
Let $L_k=\frac12\left(1-\frac1{2^k}\right)$. We take $\psi(x)=\psi_{L_k}(x)=L_k+(|x|^{\sigma/2}-1)_+$ in~\eqref{eq:newenergy},  and put $w_k(t)=(w-\psi_{L_k})_{+}(\cdot,t)$. Since $\psi_k\ge0$ we take $\ell=0$ in~\eqref{eq:lambdas}. Observe that if we start the iteration from $k=1$ we may take $\ell=1/4$.  Also,
when $w>\psi_k$, condition~\eqref{eq:primerapsi} implies $w\le\frac{1+\sqrt5}2$. We take $M=2$ in~\eqref{eq:lambdas} to simplify.

We define the quantity
\begin{equation*}\label{eq:Uk}
U_{k}=\sup_{t_{k}<t<0}\|w_{k}(t)\|_2^{2}+ \int_{t_{k}}^{0}\|(-\Delta)^{\sigma/4}w_k(t)\|_2^2\,dt,\qquad t_k=-1-\frac{1}{2^k}.
\end{equation*}
The energy estimate~\eqref{eq:newenergy} implies, for $k\ge1$, that
\begin{equation}\label{eq:Uk2}
U_k\le C2^k\,\frac{1+\Lambda_2}{\Lambda_1}\int_{t_{k-1}}^{t_k}\left(\|w_k(t)\|_2^2+
\|w_k(t)\|_1+\|\mathds{1}_{\{w_k(t)>0\}}\|_1\right)dt.
\end{equation}

Now, since $L_k=L_{k-1}+2^{-k-1}$, we have that $w_k>0$ implies $w_{k-1}>2^{-k-1}$, which in turn gives the Chebyshev type inequality
$$
\int_{\mathbb{R}^N}w_k^p\le2^{(k+1)(q-p)}\int_{\mathbb{R}^N}w_{k-1}^q
$$
for every $q>p$. Thus, for some $q>2$ to be chosen we get that \eqref{eq:Uk2} reduces to
\begin{equation*}\label{eq:Uk3}
\begin{array}{rl}
U_k&\displaystyle\le C2^k\,\frac{1+\Lambda_2}{\Lambda_1}\int_{t_{k-1}}^{t_k}\left(2^{(k+1)(q-2)}+2^{(k+1)(q-1)}+2^{(k+1)q}\right)\int_{\mathbb{R}^N}w_{k-1}^q(t)dt
\\ [4mm]
&\displaystyle\le C2^{(q+1)k}\,\frac{1+\Lambda_2}{\Lambda_1}\int_{t_{k-1}}^{t_k}\|w_{k-1}(t)\|_q^qdt.
\end{array}
\end{equation*}

To link this estimate with $U_{k-1}$, usually Hardy-Littlewood-Sobolev inequality~\eqref{HLS} is used, so $N>\sigma$ is required. The following nonlinear recurrence
\begin{equation}\label{eq:ukuk-1}
U_k\le (C_{N,\sigma,\Lambda})^k \frac{1+\Lambda_2}{\Lambda_1}\,U_{k-1}^{1+\frac{\sigma}{N}},
\end{equation}
is obtained.
Hence we are left with the case $\sigma\ge N$, which is only possible if $N=1$. The idea to deal with this range of parameters is to substitute Hardy-Littlewood-Sobolev inequality by the Nash-Gagliardo-Nirenberg inequality~\eqref{NGN}. Using first interpolation and then~\eqref{NGN}   we get, with $q=2(1+\sigma)$,
$$
\begin{array}
[c]{rl}%
\displaystyle\int_{t_{k-1}}^{t_k}\|w_{k-1}(t)\|_{q}^{q}\,dt & \displaystyle\le
\int_{t_{k-1}}^{t_k}\|w_{k-1}(t)\|_{2}^{2(\sigma-1)}\|w_{k-1}(t)\|_{\frac
{4}{2-\sigma}}^{4}\,dt\\[4mm]
& \displaystyle\le\left( \sup_{t_{k-1}<t<0}\|w_{k-1}(t)\|_{2}^{2}\right)
^{\sigma} \int_{t_{k-1}}^{t_k}\|(-\Delta)^{\sigma/4}w_{k-1}(t)\|_{2}^2\,dt\\[4mm]
& \displaystyle\le C\left( \sup_{t_{k-1}<t<0}\|w_{k-1}(t)\|_{2}^{2}+ \int_{t_{k-1}%
}^{t_k}\|(-\Delta)^{\sigma/4}w_{k-1}(t)\|_2^2\,dt\right)
^{1+\sigma}\\[6mm]
& \displaystyle\le C U_{k-1}^{1+\sigma}.
\end{array}
$$
We get again \eqref{eq:ukuk-1}. Thus, if $\frac{1+\Lambda_2}{\Lambda_1}U_0^{\sigma/N}$ is small, i.e.,
\begin{equation}\label{eq:u0small}
\int_{-2}^0\int_{\mathbb{R}^N}\big(w-(|x|^{\sigma/2}-1)_+\big)_+^2<
\varepsilon\left(\frac{\Lambda_1}{1+\Lambda_2}\right)^{1+N/\sigma},
\end{equation}
then  $U_k\to0$ as $k\to\infty$, which gives
\begin{equation}\label{w<psi}
w(x,t)\le \frac12+(|x|^{\sigma/2}-1)_+\quad\text{for }x\in\mathbb{R}^N,\;-1<t<0.
\end{equation}

The result now follows from a scaling argument. Let $(x_0,t_0)\in \Gamma_1$ be arbitrary, and define for some large $R$ the function
\begin{equation*}
  \label{eq:wR}
  w_R(x,t)=w(x_0+R^{-1}x,t_0+R^{-\sigma}t).
\end{equation*}
This function solves the equation $\partial_t\vartheta(w_R)+\mathcal{K}_Rw_R=0$, where $\mathcal{K}_R$ is the nonlocal integral operator associated to the rescaled kernel
$$
J_R(x,y)=R^{-(N+\sigma)}J(x_0+R^{-1}x,x_0+R^{-1}y).
$$
Observe that this kernel satisfies again hypothesis~\eqref{eq:kernel} with the same constant provided $R\ge1$. We now study the condition~\eqref{eq:u0small} for this function $w_R$. Observe that $1+(|x/R|^{\sigma/4}-1)_+\le |x|^{\sigma/2}-1$ for every $|x|>R$ if $R$ is large enough. Thus
$$
\begin{array}{rl}
\displaystyle\int_{-2}^0\int_{\mathbb{R}^N}\big(w_R(x,t)-(|x|^{\sigma/2}-1)_+\big)_+^2\,dxdt&\displaystyle
\le\int_{-2}^0\int_{B_R(0)}(w_R)_+^2 \\ [10pt]
&\displaystyle\le R^{N+\sigma}\int_{t_0-2/R^\sigma}^{t_0}\int_{B_1(x_0)}w_+^2 \\ [12pt]
&\displaystyle
\le R^{N+\sigma}2^{\sigma/4}|\{w>0\}\cap \Gamma_2|,
\end{array}
$$
since $B_1(x_0)\times(t_0-2R^{-\sigma},t_0)\subset\Gamma_2$ if $R^\sigma>2$, and $w\le2^{\sigma/4}$ in $B_2(x_0)$, thanks to~\eqref{eq:primerapsi}.
Choosing $c=\varepsilon R^{-(N+\sigma)}2^{-\sigma/4}$ in \eqref{eq:critical.condition}, we get from \eqref{w<psi} that $w_R<1/2$ in $\Gamma_1$, which implies $w(x_0,t_0)<1/2$.
\qed

\noindent\emph{Remark. }
As it is noted in the proof, the result also holds in terms of the constant
\begin{equation}
\label{eq:delta.vartheta2}
\overline{\delta}(\vartheta):=\frac{\inf\limits_{1/4\le s\le 2}\vartheta'(s)}{1+\vartheta(2)-\vartheta(1/4)}>0,
\end{equation}

\medskip

To proceed with the regularity proof we need to analyse what happens when the solution is neither mostly negative nor mostly positive, in the sense of Lemma~\ref{lem:first.DG.lemma}, in space-time measure. To this aim we will use De Giorgi's idea of loss of mass at intermediate levels,  obtaining a quantitative version of the fact that a function with a jump discontinuity cannot be in the energy space.

The key idea is to impose conditions on the nonlinearity guaranteeing that the equation is neither degenerate nor singular at the intermediate values. Hence we are in the linear setting studied in~\cite{Caffarelli-Chan-Vasseur}. As there, the result is written in terms of the functions
$$
\psi_\lambda(x)=((|x|-\lambda^{-4/\sigma})_+^{\sigma/4}-1)_+, \qquad \lambda\in (0,1/3),
$$
used to control the growth at infinity, and
$$
F(x)=\sup(-1,\inf(0,|x|^2-9)),
$$
used to \lq\lq localize'' the problem in the ball $B_3$. Notice that $F$
equals -1 in $B_1$, and vanishes outside $B_3$.

\begin{Lemma}
\label{lem:second.DG.lemma}
Assume $C_1\le \vartheta'(s)\le C_2$ for every $1/2\le s\le2$.
For every $\nu,\mu>0$  there exist $\gamma>0$ and $\bar\lambda\in(0,1/3)$ such that for any $\lambda\in (0,\bar\lambda)$, and any solution $w:\mathbb{R}^N\times[-3,0]\to\mathbb{R}$ to~\eqref{eq:main-vartheta} satisfying
$$
w(x,t)\le 1+\psi_\lambda(x)\quad\text{on }\mathbb{R}^N\times[-3,0],\qquad |\{w<0\}\cap(B_1\times (-3,-2))|>\mu,
$$
we have the following implication: If
$$
|\{w>1+\lambda^2F\}\cap (B_3\times(-2,0))|\ge\nu,
$$
then
$$
|\{1+F<w<1+\lambda^2F\}\cap (B_3\times(-3,0))|\ge\gamma.
$$
\end{Lemma}

The main idea in the linear case is that the  truncation function
$$
\bar\psi=1+\psi_\lambda+\lambda F
$$
satisfies an improved energy estimate. Since $\bar\psi\ge1-\lambda>1/2$, and $w\le1+\psi_\lambda=1$ in $B_3$, our assumptions on the nonlinearity $\vartheta'$ give
$$
\frac{C_1}2(w-\bar\psi)_+^{2}\le
\mathcal{B}_{\bar\psi}(w)\le\frac{C_2}2(w-\bar\psi)_+^{2},
$$
and the proof in~\cite[Lemma 4.1]{Caffarelli-Chan-Vasseur} works verbatim, using the above equivalence whenever required.

We have now all the ingredients to prove the oscillation reduction result. The growth at infinity is controlled in this case by
\[
H_{\lambda,\varepsilon}(x)=[(|x|-c(\lambda))^{\varepsilon}-1]_{+},
\]
with $\lambda>0$ small, $c(\lambda)$ large and $\varepsilon>0$.

\begin{Lemma}\label{lema-oscil}
Let $\vartheta$ be such that $\delta(\vartheta)>0$ and $\delta(\widetilde\vartheta)>0$, where $\widetilde\vartheta(s)=-\vartheta(-s)$.
Assume in addition that $C_1\le\vartheta'(s)\le C_2$ for $s\in[1/2,2]$ or $s\in[-2,-1/2]$. There exist constants $\varepsilon>0$ and $\lambda^{*}\in(0,1)$ such that if
$w$  is a solution  to~\eqref{eq:main-vartheta} that   satisfies, for $\lambda\in(0,\bar\lambda)$ small enough,
\[
|w(x,t)|\le 1+H_{\lambda,\varepsilon}(x)\qquad\mbox{for every } x\in
\mathbb{R}^{N},\;-3\le t\le 0,
\]
then
\[
\sup\limits_{\Gamma_1} w-\inf\limits_{\Gamma_1} w\le2-\lambda^{*}.
\]
\end{Lemma}

\noindent\emph{Proof. } If  $w$ (or $-w$) is subcritical at the level $0$, i.e., if $|\{w>0\}\cap \Gamma_2|\le c\delta(\vartheta)^{1+N/\sigma}$, see~\eqref{eq:critical.condition}, we are done thanks to Lemma~\ref{lem:first.DG.lemma}. Notice that $-w$ solves~\eqref{eq:main-vartheta} with $\vartheta$ replaced by $\widetilde\vartheta$. Otherwise, thanks to the hypotheses on $\vartheta'$, either $w$ or $-w$ satisfies the hypotheses of Lemma~\ref{lem:second.DG.lemma}. We assume for definiteness that it is $w$.

We consider the sequence of rescaled functions
$$
w_{k+1}=\frac{w_k-(1-\lambda^2)}{\lambda^2},\qquad w_0=w.
$$
We have that $w_k$ is a weak solution of problem~\eqref{eq:main-vartheta} with a nonlinearity $\vartheta_{k+1}$ given iteratively by
$$
\vartheta_{k+1}(s)=\frac{1}{\lambda^2}\vartheta_k(\lambda^2 s+ 1-\lambda^2),\qquad \vartheta_0=\vartheta,
$$
always with the same operator $\mathcal{K}$. We will prove that for each $k$ we can apply either Lemma~\ref{lem:first.DG.lemma} or Lemma~\ref{lem:second.DG.lemma}. Repeated application of Lemma~\ref{lem:second.DG.lemma} will give that in fact Lemma~\ref{lem:first.DG.lemma} can be applied after a finite number of steps. Hence we will be done.

The key point is that $\vartheta_{k+1}'(s)=\vartheta_k'(\lambda^2s+1-\lambda^2)$. Hence, on the one hand, since $\lambda^2s+1-\lambda^2\in[1/2,2]$ whenever $s\in[1/2,2]$, we have  $C_1\le \vartheta_k'(s)\le C_2$ for every $k$. On the other hand,
since $[1-\lambda^2,1+\lambda^2]\subset[1/2,2]$, we get $\delta(\vartheta_k)\ge \bar\delta>0$ for all $k$.

Let $\nu=c\bar\delta^{1+N/\sigma}$, with $c$ as in Lemma~\ref{lem:first.DG.lemma}. Assume by contradiction that no $w_k$ is subcritical, that is, $|\{w_k>0\}\cap \Gamma_2|> \nu$ for all $k$, so that we could never apply Lemma~\ref{lem:first.DG.lemma}. Let $\mu>0$ be such that $|\{w<0\}\cap (B_1\times(-3,-2))|\ge\mu$. By construction,
$$
|\{w_{k+1}<0\}\cap (B_1\times(-3,-2))|\ge |\{w_k<0\}\cap (B_1\times(-3,-2))| \ge\mu.
$$
We chose $\lambda$ and $\varepsilon$ small enough, so that
$\frac{H_{\varepsilon,\lambda}(x)}{\lambda^{2}}\le \psi_\lambda(x)$.
Since
$$
w_{k+1}(x,t)\le 1+\frac{w_k(x,t)}{\lambda^2},
$$
we get by induction that $w_k(x,t)\le 1+\psi_\lambda(x)$. Then, applying Lemma~\ref{lem:second.DG.lemma}
$$
\begin{array}{l}
|\{w_{k+1}>1+\lambda^2 F\}\cap \Gamma_3|\\
\qquad\qquad
=|\{w_{k+1}>1+F\}\cap \Gamma_3|-|\{1+\lambda^2F>w_{k+1}>1+F\}\cap \Gamma_3|\\
\qquad\qquad\le|\{w_k>1+\lambda^2 F\}\cap \Gamma_3|-\gamma\le |\{w>1+\lambda^2 F\}\cap \Gamma_3|-k\gamma,
\end{array}
$$
and we arrive to a contradiction if  $k\ge|\Gamma_3|/\gamma$. We conclude that
$$
w_{k_*}\le \frac12\quad\text{in }\Gamma_1 \quad\text{for some } k_*\le |\Gamma_3|/\gamma.
$$
Going back to the original variables we get that $w= 1+\lambda^{2k_*}(w_{k_*}-1)\le1-\lambda^*$, $\lambda^*=\lambda^{2|\Gamma_3|/\gamma}/2$.
\qed

This result  shows in particular that the oscillation of $w$ in $\Gamma_2$ is reduced in $\Gamma_1$ by a factor $\varpi^*=1-\lambda^*/2$. From this we get next the regularity stated in  Theorem~\ref{th:regularity} by means of scaling arguments. As in~\cite{Caffarelli-Soria-Vazquez}, 
we have to consider separately the degenerate and nondegenerate cases.

\noindent\emph{Proof of Theorem~\ref{th:regularity}. }
\textsc{Normalization. }
Let $(x_0,t_0)\in Q$ and $\tau_0=\inf\{1,t_0/3\}$. Then
$$
v_0(x,t)=\frac{w(x_0+\tau_0^{1/\sigma}x,t_0+\tau_0t)}{\|w(\cdot,0)\|_\infty}
$$
is a solution to the equation
$$
\partial_t \vartheta_0(v_0)+\mathcal{K}_0 v_0=0,
$$
in $\mathbb{R}^N\times(-3,0)$, where $\vartheta_0(s)=\frac{1}{\|w(\cdot,0)\|_\infty}\beta(\|w(\cdot,0)\|_\infty s)$, and the operator $\mathcal{K}_0$ is the nonlocal
integral operator associated to the rescaled kernel
$$
J_0(x,y)=\tau_0^{\frac{N+\sigma}\sigma}J(x_0+\tau_0^{1/\sigma}x,x_0+\tau_0^{1/\sigma}y).
$$
The function $\vartheta_0$ and the operator $\mathcal{K}_0$ satisfy the
hypotheses of Lemma~\ref{lema-oscil}.

\noindent\textsc{Modulus of continuity. } We prove that $v_0$ is continuous at $(0,0)$.
Given $R>1$, we define the
sequence of functions, for $k\ge1$,
$$
\displaystyle v_{k}(x,t)=\frac{v_0(R^{-(k+1)}x,R^{-\sigma (k+1)}t)-\mu_{k}}{\varpi_{k}},
$$
where $\varpi_k$ and $\mu_k$ are respectively the semi-oscillation and a certain mean of $v_0$ in the parabolic cylinder $Q_k=\Gamma_{R^{-k}}$,
$$
\displaystyle\varpi_k=\frac{\sup_{Q_k} v_0-\inf_{Q_k} v_0}{2},\qquad
\mu_k=\frac{\sup_{Q_k} v_0+\inf_{Q_k} v_0}{2}.
$$
They satisfy the equation
$$
\partial_t\vartheta_k(v_k)+\mathcal{K}_kv_k=0, \qquad\vartheta_{k}(s)=\dfrac{\vartheta_0\left(\varpi_{k} s+\mu_{k}\right)}{\varpi_{k}},
$$
where the operator $\mathcal{K}_k$ has associated kernel
$$
J_{k}(x,y)=R^{-(N+\sigma)(k+1)}J_0(R^{-(k+1)}x,R^{-(k+1)}y),
$$
that satisfies again~\eqref{eq:kernel}.
Assuming by contradiction that $\varpi_{k}\ge \varsigma>0$, we have that the function $\vartheta_k$  satisfies the hypotheses of
Lemma~\ref{lema-oscil}, since $\varpi_ks+\mu_k\ge \varsigma/2$ for $s\ge1/2$ if $\mu_k\ge0$ and $\varpi_ks+\mu_k\le -\varsigma/2$ for $s\le-1/2$ if $\mu_k\le0$.

On the other hand, $|v_k|\le 1\le 1+H_{\lambda,\varepsilon}(x)$ for $|x|\le R$, applying
Lemma~\ref{lema-oscil} by induction to $v_{k-1}$, since it can be applied to $v_0$. Also, if we take $R>1$ large enough so that
$H_{\lambda,\varepsilon}(R)\ge\frac{2-\varsigma}{\varsigma}$, we get $|v_k(x,t)|\le 1+H_{\lambda,\varepsilon}(x)$ if $|x|\ge R$.  Hence, applying
Lemma~\ref{lema-oscil} we conclude that
$\varpi_k\le(1-\lambda^*/2)^k$, a contradiction. Therefore we have a modulus of continuity.

\noindent\textsc{H\"older regularity at nondegeneracy points. } We assume $w(x_0,t_0)>0$, the case $w(x_0,t_0)<0$ being similar.
We define now iteratively the
sequence of functions
$$
\displaystyle v_{k+1}(x,t)=\frac{v_k(R^{-1}x,R^{-\sigma}t)-\mu_{k}^*}{\varpi^*},
$$
where
$$
\displaystyle\varpi^*=1-\lambda^*/2,\qquad
\mu_k^*=\frac{\sup_{Q_1} v_k+\inf_{Q_1} v_k}{2}.
$$
Observe that the recurrence relation can be written explicitly,
$$
\displaystyle v_{k}(x,t)=\frac{v_0(R^{-(k+1)}x,R^{-\sigma (k+1)}t)-\nu_{k}}{(\varpi^*)^{k}},
$$
where $\nu_k=\sum\limits_{j=1}^k\mu_j^*(\varpi^*)^k$. Also, $\mu_k=(\varpi^*)^k\mu_{k+1}^*+\nu_k$,
so that, since $\mu_k\to \frac{w(x_0,t_0)}{\|w(\cdot,0)\|_\infty}>0$, then $\nu_k\to a>0$.

The functions $v_k$  satisfy $|v_k|\le 1$ in $\Gamma_R$, and the equation
$$
\partial_t\vartheta_k(v_k)+\mathcal{K}_kv_k=0,
$$
where the new nonlinearity is
$$
\vartheta_{k}(s)=\dfrac{\vartheta_0\left((\varpi^*)^ks+\nu_k\right)}{(\varpi^*)^k}
$$
and the operator $\mathcal{K}_k$ is as before.
The function $\vartheta_k$ and the operator $\mathcal{K}_k$ satisfy once more the hypotheses of
Lemma~\ref{lema-oscil}.

On the other hand, if we take $R>1$ large enough so that
$$
H_{\lambda,\varepsilon}(x/R)\le \frac{\varpi^*}2H_{\lambda,\varepsilon}(x), \quad H_{\lambda,\varepsilon}(x)\ge\frac{2(2-\varpi^*)}{\varpi_*}\quad\text{for }|x|\ge R,
$$
then
$|v_k(x,t)|\le 1+H_{\lambda,\varepsilon}(x)$ if $|x|\ge R$. We conclude, applying
Lemma~\ref{lema-oscil}, an oscillation estimate of order
$(\varpi^*)^k$ for $w$ in $Q_k$. This gives  H\"older regularity at points where the equation is nondegenerate.

\noindent\textsc{H\"older regularity at degeneracy points. }
Let now $w(x_0,t_0)=0$.
Here we consider the
sequence of functions defined by means of a recurrence that takes into account the nonlinearity, and the possible singularity of $\beta'$ at zero:
$$
\displaystyle v_{k+1}(x,t)=\frac{v_k(R^{-1}x,\gamma R^{-\sigma}t)-\mu_{k}^*}{\varpi^*},\qquad \gamma=\frac{\vartheta_0(\varpi^*)}{\varpi^*},
$$
with $\mu_k^*$ and $\varpi^*$ as before.
The rescaled nonlinearity turns to be
$$
\vartheta_{k}(s)=\dfrac{\vartheta_0\left((\varpi^*)^ks+\nu_k\right)}{\vartheta_0((\varpi^*)^k)}.
$$
We observe that
\begin{equation}\label{eq:uncociente}
\frac{|\nu_k|}{(\varpi^*)^k}\le \frac{|\mu_k|}{(\varpi^*)^k}+|\mu_{k+1}^*|\le C.
\end{equation}

The conditions of Lemma~\ref{lema-oscil} are fulfilled as long as, for every $k\ge1$,
$$
0<C_1\le \frac{(\varpi^*)^k\vartheta_0'((\varpi^*)^ks+\nu_k)}{\vartheta_0((\varpi^*)^k)}\le C_2\quad\text{for every }s\in(1/2,2),
$$
and
$$
\frac{\frac{(\varpi^*)^k}{\vartheta_0((\varpi^*)^k)}
\inf\limits_{[\nu_k,2(\varpi^*)^k+\nu_k]}\vartheta_0'}{1+\frac{\vartheta_0(2(\varpi^*)^k+\nu_k)
-\vartheta_0(\nu_k)}{\vartheta_0((\varpi^*)^k)}}\ge \ell>0.
$$
They hold from condition~\eqref{eq:cond.phi01} using \eqref{eq:uncociente}.
We conclude as before.
\qed


\section{Asymptotic behaviour. Proof of Theorem~\ref{th:asymptotic.behaviour}}

\label{sect-asymptotic-behaviour} \setcounter{equation}{0}

We devote this section to study the large time behaviour of bounded solutions
to~\eqref{eq:main}. Notice that, since the solution is bounded, we may assume (modifying the nonlinearity for large values of $u$, if required), that there exist constants $0<c\le C<\infty$ such that
 \begin{equation}
\label{eq:behaviour.phi'.u.bounded}
c |u|^{m-1}\le \varphi'(u)\le C |u|^{m-1}\qquad\text{ for all }u\in\mathbb{R},
\end{equation}
since this is true for $u\approx0$. We may assume also, for simplicity, by a simple rescaling, that $a=m$ and $\mu=\mu_{N,\sigma}$ in~\eqref{eq:cond.J},~\eqref{eq:cond.phi.prima}. On the other hand, by Theorem~\ref{th:regularity}, we may assume that $u$ is  H\"older continuous for $t\ge0$.

We use the nowadays classical method of scalings. Let us
consider the sequence of functions
\[
u_{k}(x,t)=k^{\alpha}u(k^{\beta} x,k t),\quad k>0,
\]
with $\alpha$ and $\beta$ as in~\eqref{eq:Barenblatt}. Notice that the scaling preserves mass. More precisely,
\begin{equation*}
\label{eq:mass.is.preserved}
\int_{\mathbb{R}^N}u_k(x,t)\,dx=\int_{\mathbb{R}^N}u(x,t)\,dx.
\end{equation*}
It is trivial to check that $u_{k}$ satisfies
\[
\partial_{t}u_{k}+{\mathcal{L}}_{k}\varphi_{k}(u_{k})=0,\qquad u_{k}%
(x,0)=k^{\alpha}u_{0}(k^{\beta}x),
\]
where $\varphi_{k}(s)=k^{m\alpha}\varphi(s/k^{\alpha})$, and  $\mathcal{L}_{k}$ has associated kernel
$J_{k}(z)=k^{\beta(N+\sigma)} \widetilde J (k^{\beta}z)$.

We first observe that the operators of the family $\mathcal{L}_{k}$ and the functions of the family $\varphi_{k}$ satisfy the hypotheses of Theorems~\ref{th:regularity} and~\ref{th:smoothing}. On the other hand, the assumptions~\eqref{eq:cond.J} and~\eqref{eq:cond.phi.prima} give
\begin{align}
\label{eq:cond.Jk}
&\displaystyle\lim_{k\to\infty}J_k(z)=J_\infty(z):=\mu_{N,\sigma}|z|^{-N-\sigma}\quad&\text{uniformly for }
|z|\ge K>0;
\\
\label{eq:cond.phik}
&\displaystyle \lim_{k\to\infty}\varphi_k(s)=\varphi_\infty(s):=|s|^{m-1}s\quad&\text{uniformly for } |s|\le K<\infty,
\\
\label{eq:cond.phik.prima}
&\displaystyle \lim_{k\to\infty}\varphi'_k(s)=\varphi_\infty'(s)\quad&\text{uniformly for }|s|\le K<\infty;
\end{align}
for some $m\ge1$. Moreover, from the lower bound in~\eqref{eq:behaviour.phi'.u.bounded},
\begin{equation}
\label{eq:cond.para.holder.uniforme}\varphi_{k}^{\prime}(s)\ge C |s|^{m-1}\quad\text{for some
constant }C>0\text{ independent of }k.
\end{equation}
Hence, since $\|u_{0,k}\|_{L^{1}(\mathbb{R}^{N})}=\|u_{0}
\|_{L^{1}(\mathbb{R}^{N})}$, the smoothing effect
\eqref{eq:L-inf-p3} tells us that
\begin{equation}
\label{eq:uniform.bound}
\|u_{k}(\cdot,t)\|_{L^{\infty}(\mathbb{R}^{N})}\le
c_{1}(\nu),\qquad t\ge\nu>0.
\end{equation}

Thanks to \eqref{eq:cond.para.holder.uniforme} and~\eqref{eq:uniform.bound}, we may apply
Theorem~\ref{th:regularity}, and we obtain
that the family
$\{u_{k}\}$ is uniformly H\"{o}lder continuous for $t\ge\nu'>\nu$. Then, applying
Ascoli-Arzela's lemma, we have that there exists a subsequence converging uniformly on
compact subsets of $Q$ to some function $v\in C^{\alpha
}(\mathbb{R}^N\times[\nu',\infty))$. Again by translation in time we may assume $\nu'=0$.
For a convergent subsequence $\{u_{k}\}$, since it is uniformly bounded,~\eqref{eq:cond.phik} yields
$\varphi_{k}(u_{k})\to |v|^{m-1}v$ uniformly in compact subsets of $Q$. If we are able to identify the limit $v$ as $B_M$, the result will follow by the classical procedure of taking $t=1$ and $k=t$. This identification is our next objective.

Given $\tau>0$, we consider the
translate in time of the Barenblatt solution to the fractional porous medium equation with mass $M=\int_{\mathbb{R}^N}u_0$,
\[
B_{M,\tau}(x,t)=B_M(x,t+\tau).
\]
Notice that $B_{M,\tau}\to B_{M}$  as $\tau\to0$ weakly in $L^{q}(Q_T)$, $Q_T=\mathbb{R}^N\times(0,T)$, $T>0$,
$q\in(1,m+\sigma/N)$. Therefore, to show that $v=B_M$ it is enough to prove that, given any
$F\in C^{\infty}_{\mathrm{c}}(Q)$,  for all $\varepsilon>0$ there exists a value
$\tau_{\varepsilon}>0$,  and for each $\tau$ a constant $K=K(\varepsilon,\tau)$ such that
\[
\left| \iint_{Q}(u_k-B_{M,\tau})F \right|
\le\varepsilon\quad\text{for } \tau\le\tau_{\varepsilon},\ k\ge K(\varepsilon,\tau).
\]

Let $f$ be an admissible test function for our rescaled nonlinear problems.
Then, for all $k>1$ and $\tau>0$ we have
\[%
\begin{array}
[c]{l}%
\displaystyle\iint_{Q}\left( (u_{k}-B_{M,\tau})\partial_{t} f- \varphi_{k}(u_{k}) \mathcal{L}_k f+B_{M,\tau}^m
(-\Delta)^{\sigma/2}f\right) \\[4pt]%
\displaystyle\qquad=-\int_{\mathbb{R}^{N}}\left( u_{k}(\cdot,0)-B_{M,\tau
}(\cdot,0)\right) f(\cdot,0).
\end{array}
\]
This may be rewritten as
\[%
\begin{array}
[c]{l}%
\displaystyle\iint_{Q} (u_{k}-B_{M,\tau})\left( \partial_{t} f-\left(a_{k,\tau}+\frac1n\right)
(-\Delta)^{\sigma/2}f\right) \\[10pt]%
\displaystyle\qquad=-\underbrace{\int_{\mathbb{R}^{N}}\left( u_{k}%
(\cdot,0)-B_{M,\tau}(\cdot,0)\right) f(\cdot,0)}_{I_{1}}\\[10pt]
\displaystyle\qquad\quad
+\underbrace{\iint_{Q}(-\Delta)^{\sigma/4}\left(\varphi_{k}(B_{M,\tau})-\varphi_\infty(B_{M,\tau})\right)
(-\Delta)^{\sigma/4}f}_{I_{2}}
\\[10pt]
\displaystyle\qquad\quad+\underbrace{\iint_{Q}\varphi_{k}(u_k)
\left(\mathcal{L}_k-(-\Delta)^{\sigma/2}\right)f}_{I_{3}}+\underbrace{\frac1n
\iint_{Q} (u_{k}-B_{M,\tau})(-\Delta)^{\sigma/2}f}_{I_{4}},
\end{array}
\]
where
\[
a_{k,\tau}=
\begin{cases}
\displaystyle\frac{\varphi_{k}(u_{k})-\varphi_{k}(B_{M,\tau})}{u_{k}-B_{M,\tau}} &
\text{if }u_{k}\ne B_{M,\tau},\\
\varphi^{\prime}_{k}(B_{M,\tau}) & \text{if }u_{k}=B_{M,\tau}.
\end{cases}
\]
Thanks to Theorem~\ref{th:regularity} and the smoothness assumptions on $\varphi$, and using the upper bound in~\eqref{eq:behaviour.phi'.u.bounded}, we know that $a_{k,\tau}\in
C^{\alpha}(Q)\cap L^{\infty}(Q)$.

Let now $f=f_{n,k,\tau}$ be a classical solution to
\[
\partial_{t} f-\left(a_{k,\tau}+\frac1n\right)(-\Delta)^{\sigma/2}f=F,\qquad f(x,T)=0.
\]
It is known that such a solution exists; see Appendix. Moreover, if extended by zero for $t>T$, it is an admissible test function. We now proceed to prove that $I_i\to0$, $i=1,\dots,4$ for this particular choice of $f$. The limit is taken as $n\to\infty$, then $k\to\infty$ and finally $\tau\to0$.

In order to estimate $I_{1}$, we use that $f$ is H\"{o}lder continuous and bounded  (uniformly in $n$, $k$ and $\tau$; see Theorem~\ref{th:cauchy}) at $t=0$, and
that $u_{k}(\cdot,0)$ and $B_{M,\tau}$ have the same integral, to obtain
$$
\begin{array}
{l}
|I_{1}|  =  \displaystyle\left|  \int_{\mathbb{R}^{N}}\left( u_{k}
(x,0)-B_{M,\tau}(x,0)\right) (f(x,0)-f(0,0))\,dx\right| \\[10pt]
\le \displaystyle C \left(R^{\alpha}\int_{|x|\le R}|u_{k}(x,0)-B_{M,\tau}(x,0)|\,dx
+\int_{|x|\ge k^{\beta}R}u_{0}(x)\,dx+ \int_{|x|\ge R}
B_{M,\tau}(x,0)\,dx\right)\\[10pt]
\le \displaystyle C\left(R^{\alpha}+\int_{|x|\ge k^{\beta}R}
u_{0}(x)\,dx+  \int_{|x|\ge R} B_{M,\tau}(x,0)\,dx\right).
\end{array}
$$
It is easily seen that we can make $|I_{1}|<\varepsilon$ just taking first $R$
small enough, and then $k$ big and $\tau$ small.

As for $I_{2}$, we have
\[
\begin{array}{l}
\displaystyle|I_{2}|\le\|(-\Delta)^{\sigma/4}\left(\varphi_{k}(B_{M,\tau})-
\varphi_\infty(B_{M,\tau})\right)\|_{2}\|(-\Delta)^{\sigma/4}f\|_{2}\\[10pt]
\displaystyle\qquad\le \sup_{0\le\theta\le\|B_\tau\|_\infty}\|\varphi_{k}^{\prime}(\theta)-\varphi'_\infty(\theta)\|_{\infty}
\|(-\Delta)^{\sigma/4}B_{M,\tau}\|_{2}\|(-\Delta)^{\sigma/4}f\|_{2}\to 0
\end{array}
\]
as $k\to\infty$ for all $\tau\in(0,\tau_\varepsilon)$ fixed,
thanks to \eqref{eq:cond.phik.prima} and the uniform estimate for $\| f\|_{\dot H^{\sigma/2}}$ in~\eqref{eq:prop-adjoint}.

The estimate for $I_3$ will follow from condition~\eqref{eq:cond.Jk} and the regularity of $f$. In fact, since the family $\varphi_k(u_k)$ is uniformly bounded in $L^1(Q_T)$, it is enough to estimate $\left(\mathcal{L}_k-(-\Delta)^{\sigma/2}\right)f$ in $L^\infty(Q_T)$.   We recall that since $J$ satisfies~\eqref{eq:J.symmetric},  we can use expression~\eqref{operatorL2} both for $\mathcal{L}_k$ and $(-\Delta)^{\sigma/2}$, thus getting
$$
\begin{array}{l}
\displaystyle\left|\left(\mathcal{L}_k-(-\Delta)^{\sigma/2}\right)f\right(x,t)|
\\ [4mm]\qquad\displaystyle\le\frac12\int_{\mathbb{R}^N}|f(x+y,t)+f(x-y,t)-2f(x,t)||J_k(y)-J_\infty(y)|\,dy \\ [10pt]
\qquad\displaystyle\le C\int_{|y|\le R}\frac{|y|^{\sigma+\varepsilon}}{|y|^{N+\sigma}}\,dy+
C\int_{|y|\ge R}|J_k(y)-J_\infty(y)|\,dy\to0
\end{array}
$$
uniformly in $Q_T$.

Finally using  estimate~\eqref{eq:prop-adjoint} we obtain

$$
\begin{array}{rcl}
|I_{4}|&\le&\displaystyle\frac1n\left( \iint_{Q} (u_{k}-B_{M,\tau})^{2}\frac
{1}{a_{k,\tau}+\frac1n}\right) ^{1/2}\left( \iint_{Q}
\left(a_{k,\tau}+\frac1n\right)((-\Delta)^{\sigma/2}f)^{2}\right) ^{1/2}\\
&\le& C/n^{1/2}.
\end{array}
$$

\section*{Appendix. Parametrix method}\label{sect-appendix}
\setcounter{equation}{0}
\newcommand{\oy}{\overline Y}
\renewcommand{\d}{\,\mathrm{d}}
\renewcommand{\theequation}{A.\arabic{equation}}
\renewcommand{\thesection}{A}

We consider the nonlocal problem in non-divergence form
\begin{equation}
\label{eq:problem.nondivergence.form}
\partial_t f+a(-\Delta)^{\sigma/2} f=F\quad\text{in }\,\mathbb{R}^N\times\mathbb{R}_+,\qquad f(\cdot,0)=f_0\in C(\mathbb{R}^N),
\end{equation}
where the coefficient $a=a(x,t)$ is H\"older continuous and  satisfies the  \lq ellipticity' condition
$$
0<\lambda_1\le a(x,t)\le \lambda_2<\infty.
$$
We moreover assume that $F=F(x,t)$ is also H\"older continuous.
Our aim is to prove that this problem is well posed in the space
$$
L_\sigma=\{g \text{ measurable } :\, \int_{\mathbb{R}^N}\frac{g(x)}{1+|x|^{N+\sigma}}\,dx<\infty\}.
$$
To this purpose we have to assume  $f_0\in L_\sigma$ and $F(\cdot,t)\in L_\sigma$ uniformly in $t\in[0,T]$.

\begin{Theorem}\label{th:cauchy}
There is a unique classical solution $f\in C(Q_T)$ to Problem \eqref{eq:problem.nondivergence.form} in $Q_T$. It also belongs to some space $C^\alpha(\mathbb{R}^N\times(\tau,T))$ for every $\tau>0$. Moreover, if
$f_0\in L^\infty(\mathbb{R}^N)\cap \dot H^{\sigma/2}(\mathbb{R}^N)$, $F\in L^\infty(Q_T)\cap  L^\infty([0,T];\dot H^{\sigma/2}(\mathbb{R}^N))$,  then the solution satisfies
\begin{equation}\label{eq:prop-adjoint}
\begin{array}{l}
\|f\|_{L^\infty(Q_{T})}\le C, \\ [4pt] 
\displaystyle\sup\limits_{0<t<T}\int_{\mathbb{R}^N}\left|(-\Delta)^{\sigma/4}f(\cdot,t)\right|^2\le  C, \\
\displaystyle\sup\limits_{0<t<T}\int_{Q_t}a\left|(-\Delta)^{\sigma/2}f\right|^2\le  C,
\end{array}
\end{equation}
where the constant $C$ depends only on  $\|f_0\|_\infty$, $\|f_0\|_{\dot H^{\sigma/2}}$, $\|F\|_\infty$ and $\sup\limits_{0<t<T}\|F(\cdot,t)\|_{\dot H^{\sigma/2}}$.
\end{Theorem}

The solution will be given by means of the representation formula
\begin{equation}
\label{sol-linear}
f(x,t)=\int_{\mathbb{R}^N}\Gamma(x,t,\xi,0)f_0(\xi)\,d\xi +
\int_0^t\int_{\mathbb{R}^N}\Gamma(x,t,\xi,\tau)F(\xi,\tau)\,d\xi
d\tau,
\end{equation}
where $\Gamma$ is the fundamental solution to the problem. It is then clear that the regularity of $f$ for $t>0$ is inherited from the regularity of $\Gamma$. Even more, the regularity of $\Gamma$ is determined by the regularity of the coefficient $a$. At this respect, assuming only local regularity, plus global boundedness, of $a$ is enough to get local regularity of $f$.

Hence, our first step is to construct $\Gamma$. This will be done adapting the parametrix method of E.~E.~Levi~\cite{Levi} to the case of the nonlocal operator $L=\partial_t+a(-\Delta)^{\sigma/2}$. The fundamental solution  has already been constructed in a very recent paper by Chen and Zhang~\cite{Chen-Zhang} that has just come to our knowledge. Nevertheless, we have decided to keep our proof since it is simpler, thanks to the use of a certain quasimetric adapted to the problem. In addition, it shows clearly the local character of the regularity result, which is in fact needed in the application to the large time behaviour of solutions to~\eqref{eq:main} given in Section~\ref{sect-asymptotic-behaviour}.

\begin{Theorem}
\label{th:parametrix}
There exists a function $\Gamma\in C(\mathbb{R}^{2N}\times\{\varepsilon \le t-\tau\le t_0\})$ satisfying $\partial_t \Gamma(\cdot,\cdot,\xi,\tau)$, $(-\Delta)^{\sigma/2}\Gamma(\cdot,\cdot,\xi,\tau)\in
C^\beta(\mathbb{R}^N\times(\tau+\varepsilon,T))$ for every $\varepsilon>0$, $T>0$ and some $\beta\in(0,1)$, for every fixed $\xi\in\mathbb{R}^{N}$, $\tau>0$, and solving
\begin{equation*}
\label{eq:fundam-sol}
\left\{
\begin{array}{ll}
\partial_t\Gamma+a(-\Delta)^{\sigma/2} \Gamma=0, & \qquad  x\in\mathbb{R}^N,\; \tau<t<T,
\\ [4mm]
\Gamma(x,\tau,\xi,\tau) = \delta(x-\xi), & \qquad x\in\mathbb{R}^N.%
\end{array}
\right.
\end{equation*}
\end{Theorem}
\noindent{\it Proof. } We construct the fundamental solution of the
operator $L=\partial_t+a(-\Delta)^{\sigma/2}$  in terms of the fundamental solution of
the operator with frozen coefficient $L_0=\partial_t+\overline
a(-\Delta)^{\sigma/2}$, $\overline a=a(\xi,\tau)$. For that purpose let
$P(x,t)$
be the Poisson kernel, which solves the fractional heat equation
$$
\partial_tP+(-\Delta)^{\sigma/2} P=0, \qquad  x\in\mathbb{R}^N,\; t>0,
$$
with $\delta(x)$ as initial value, and $P(x,t)=0$ for $t<0$, and define the function
\begin{equation}\label{Z(x,t)}
Z(x,t,\xi,\tau)=P(x-\xi,a(\xi,\tau)(t-\tau)).
\end{equation}
Then $Z(\cdot,\cdot,\xi,\tau)$ solves the fractional heat equation
with constant coefficient
$$
L_0Z=\partial_tZ+\overline a(-\Delta)^{\sigma/2} Z=0, \qquad x\in\mathbb{R}^N,\;
t>0.
$$
This function $Z$, called the \emph{parametrix}, will be the principal part of the desired function $\Gamma$. We look for $\Gamma$ in the form
\begin{equation}\label{Gamma}
\Gamma(x,t,\xi,\tau)=Z(x,t,\xi,\tau)+\int_\tau^t\int_{\mathbb{R}^N}
Z(x,t,\lambda,\eta)\Phi(\lambda,\eta,\xi,\tau)\,d\lambda d\eta,
\end{equation}
where we must construct $\Phi$ in order to have $L\Gamma=0$. The initial value is, formally,
$$
\Gamma(x,\tau,\xi,\tau)=Z(x,\tau,\xi,\tau)=\delta(x-\xi).
$$
Let now $\psi_0=-LZ=(\overline a-a)(-\Delta)^{\sigma/2} Z$. If $\Phi$ is a fixed point of the functional
$$
\mathcal{T}(\psi)(x,t,\xi,\tau)=\psi_0(x,t,\xi,\tau)+
\int_\tau^t\int_{\mathbb{R}^N}\psi_0(x,t,\lambda,\eta)\psi(\lambda,\eta,\xi,\tau)\,d\lambda d\eta,
$$
then the formal application of the operator $L$ to $\Gamma$
gives
\begin{equation}\label{formalLGamma}
L\Gamma=LZ+\Phi+\iint LZ\Phi= -\psi_0+\Phi-\iint
\psi_0\Phi=\Phi-\mathcal{T}(\Phi)=0.
\end{equation}

To construct $\Phi$ we define the
sequence $\{\psi_k\}$ for $k\ge1$ by the recurrence relation
$$
\psi_{k+1}(x,t,\xi,\tau)=\int_\tau^t\int_{\mathbb{R}^N}
\psi_0(x,t,\lambda,\eta)\psi_k(\lambda,\eta,\xi,\tau)\,d\lambda d\eta.
$$
The function $\Phi=\sum_{k=0}^\infty\psi_k$ is clearly a fixed point of $\mathcal{T}$, and thus defines $\Gamma$.
This is the standard construction performed in \cite{Levi} of the fundamental solution, described for instance in \cite{Friedman}. What remains is the justification of the calculations, i.e., the
convergence of the integrals involved as well as the convergence of
the sum, and this is specially delicate in our situation.
We first show that the sequence $\{\psi_k\}$ is well defined and that the
sum  is convergent. To that
purpose we  will use the notation $Y=(x,t),\,\overline Y=(\xi,\tau)\in Q$, as well as  the quasimetric $|Y-\oy|_\sigma$, introduced in \cite{VPQR},
\begin{equation}
\label{sigma-distance}
|Y|_\sigma:=\Big(|x|^2+|t|^{2/\sigma}\Big)^{1/2}.
\end{equation}
We also consider the corresponding H\"older space $C^\alpha_\sigma(Q)$.
In terms of the joint variables $Y$ and $\oy$, the function $Z$ in \eqref{Z(x,t)} can be written as
\begin{equation*}\label{Z(Y)}
Z(Y,\oy)=P(\vartheta(a(\oy))(Y-\oy)),
\end{equation*}
where $\vartheta:\mathbb{R}\to \mathcal{M}^{N+1}$ is given by
$$
\vartheta(s)=\left(\begin{array}{cc}
I & 0 \\
0 & s \\
\end{array}
\right),
$$
$I\in\mathcal{M}^{N}$ being the identity matrix. With this notation the fundamental solution is
\begin{equation}\label{Gamma(Y)}
\Gamma(Y,\oy)=Z(Y,\oy)+\int_QZ(Y,Y')\Phi(Y',\oy)\,dY'.
\end{equation}
The behaviour of the coefficient function $a$ between two positive constants makes the coefficient matrix $\vartheta(a(\oy))$ play no role in the estimates of $Z$ needed to study $\Gamma$. Nevertheless, writing $Z(Y,\oy)=W(Y,\oy,a(\oy))$, where $W(Y,\oy,s)=P(\vartheta(s)(Y-\oy))$, allows to estimate the behaviour when moving the variable $s$. This will be of use later on.
Now, using \eqref{sigma-distance} we can deduce from Proposition 2.1 in \cite{VPQR} the estimate
\begin{equation}\label{eq:est-lambdaZ}
|(-\Delta)^{\sigma/2} Z(Y,\oy)|\le
\dfrac{c}{|Y-\overline Y|_\sigma^{N+\sigma}}.
\end{equation}

This is not enough for  $(-\Delta)^{\sigma/2} Z$  to be integrable near $Y=\oy$. Nevertheless, the
regularity of the coefficient $(a-\overline a)$ solves this problem. The condition $a\in C^\alpha(Q_T)$ means $a\in C^{\alpha'}_\sigma(Q_T)$, for some $\alpha'=\alpha'(\alpha,\sigma)$, though we  still use the same letter $\alpha$.  Thus we have
\begin{equation*}\label{eq:est-LZ}
|\psi_0(Y,\oy)|\le
\dfrac{c}{|Y-\overline Y|_\sigma^{N+\sigma-\alpha}}
\end{equation*}
for $|Y-\oy|_\sigma$ small. For the integrability for $|Y-\oy|_\sigma$ large, since the time interval is bounded, we use the estimate
\begin{equation*}\label{eq:est-LZ2}
|\psi_0(Y,\oy)|\le
\dfrac{c}{|x-\xi|^{N+\sigma}}.
\end{equation*}
To estimate now $\psi_k$ we observe that, for $|x-\xi|<1$ we have
$$
\begin{array}{rl}
\displaystyle |\psi_1(Y,\overline Y)| &\displaystyle\le c_1\int_\tau^t\int_{|x-\lambda|<2}\big|
\psi_0(Y,Y')\big|\,\big|\psi_0(Y',\overline Y)\big|\,dY'\\
[4mm] &\displaystyle+c_1\int_\tau^t\int_{|x-\lambda|>2}\big|
\psi_0(Y,Y')\big|\,\big|\psi_0(Y',\overline Y)\big|\,dY' \\
[4mm]
&\displaystyle\le c_1\int_{|Y-Y'|<2}
\frac1{|Y-Y'|_\sigma^{N+\sigma-\alpha}}\frac1{|Y'-\overline Y|_\sigma^{N+\sigma-\alpha}}\,dY' \\
[4mm] &\displaystyle +c_2\int_{|x-\lambda|>2}\frac1{|x-\lambda|^{N+\sigma}}\frac1{|\lambda-\xi|^{N+\sigma}}\,d\lambda\\
[4mm] &\displaystyle \le\frac{c_1}{|Y-\overline Y|_\sigma^{N+\sigma-2\alpha}}+c_2.
\end{array}$$
We have used Lemma 1.2 of \cite{Friedman}. On the other hand, for $|x-\xi|>1$ we get
$$
\begin{array}{rl}
\displaystyle |\psi_1(Y,\overline Y)|&\le \displaystyle c_1\int_\tau^t\int_{|x-\lambda|<1/2}\big|
\psi_0(Y,Y')\big|\,\big|\psi_0(Y',\overline Y)\big|\,dY'\\
[4mm] &\displaystyle+c_2\int_\tau^t\int_{|x-\lambda|>1/2}\big|
\psi_0(Y,Y')\big|\,\big|\psi_0(Y',\overline Y)\big|\,dY'\le \frac{c}{|x-\xi|^{N+\sigma}}.
\end{array}$$
Therefore,
$$
|\psi_k(Y,\overline Y)|\le \frac{c_k}{|Y-\overline Y|_\sigma^{N+\sigma-(k+1)\alpha}},\qquad
|\psi_k(Y,\overline Y)|\le\frac{c}{|x-\xi|^{N+\sigma}}\quad\text{for $|x-\xi|$ large}.
$$
This means that there is a finite $k_0$ such that $\psi_k$ possesses no singularity at the origin for $k\ge
k_0$. Even more, it is easy to check that for large $k$ we have the
estimate $|\psi_k|\le \dfrac{c^k}{k!}$. This means that the sum
$\Phi=\sum_{k=0}^\infty\psi_k$ is absolutely convergent in compact subsets of
$Q_T$, it is a fixed point of the functional
$\mathcal{T}$, and satisfies the estimates
\begin{equation}
\label{eq:est-chi}
|\Phi(Y,\overline Y)|\le
\frac{c}{|Y-\overline Y|_\sigma^{N+\sigma-\alpha}},
\qquad
|\Phi(Y,\overline Y)|\le\frac{c}{|x-\xi|^{N+\sigma}}
\quad\text{for $|x-\xi|$ large}.
\end{equation}
Also, it is easy to see that
$\Phi$ is continuous in $Y$ uniformly in $\oy$ provided $t-\tau\ge c>0$. In the same way $\Phi$ is continuous in $\oy$ uniformly in $Y$ provided $t-\tau\ge c>0$. To prove that it is
H\"older continuous, we use the formula
\begin{equation}
\label{eq:fixed-chi}
\Phi(Y,\overline Y)=\psi_0(Y,\overline Y)+
\int_{Q_T}\psi_0(Y,Y')\Phi(Y',\overline Y)\,dY'.
\end{equation}
\begin{Lemma}\label{lem-holder} In the above hypotheses,
\begin{equation*}\label{eq:holder}
|\Phi(Y,\overline Y)-\Phi(\widetilde{Y},\overline Y)|\le c|Y-\widetilde{Y}|_\sigma^\alpha
\end{equation*}
for every $Y,\,\widetilde{Y},\,\oy\in Q_T$ such that $t-\tau\ge c>0$, $|Y-\widetilde{Y}|_\sigma\le c/2$.
\end{Lemma}
\noindent{\it Proof. }
The first term in \eqref{eq:fixed-chi} is H\"older continuous, since
\begin{equation}\label{I1}
\begin{array}{rl}
|\psi_0(Y,\overline Y)-\psi_0(\widetilde{Y},\overline Y)|&\le|a(Y)-a(\widetilde{Y})||(-\Delta)^{\sigma/2}
Z(Y,\oy)|
\\ [3mm]
&\quad+ |a(\widetilde{Y})-a(\oy)||(-\Delta)^{\sigma/2} Z(Y,\oy)-(-\Delta)^{\sigma/2}
Z(\widetilde{Y},\oy)|
\\ [3mm]
&\le
c\Big(\dfrac{|Y-\widetilde{Y}|^\alpha_\sigma}{|Y-\oy|_\sigma^{N+\sigma}}+
\dfrac{|Y-\widetilde{Y}|_\sigma}{|\theta-\oy|_\sigma^{N+\sigma+1}}\Big),
\end{array}
\end{equation}
by using again Proposition 2.1 in \cite{VPQR}, where $\theta$ is some intermediate point between $Y$ and $\widetilde{Y}$. Thus, since the condition $t-\tau\ge c$ implies $|Y-\oy|_\sigma\ge c$ and $|\theta-\oy|_\sigma\ge c/2$, we get
$$
|\psi_0(Y,\overline Y)-\psi_0(\widetilde{Y},\overline Y)|\le c|Y-\widetilde{Y}|^\alpha.
$$
As to the second term in \eqref{eq:fixed-chi}, we
combine \eqref{I1} with \eqref{eq:est-chi} to get the desired result. \qed

Now, in order to study the second term in \eqref{Gamma(Y)}, called the \emph{volume potential} of $\Phi$, we consider the volume potential for any given function $h:Q_T\to\mathbb{R}$,
\begin{equation*}
\label{eq:volpot}
V(h)(Y)=\int_{Q_T}Z(Y,Y')h(Y')\,dY'.
\end{equation*}
\begin{Lemma}\label{lem-volumepot}
If $h\in C^\beta_\sigma(Q_T)\cap L_\sigma(Q_T)$ for some $\beta\in(0,1)$, then
\begin{align}
\label{eq:Lambda V}
&(-\Delta)^{\sigma/2} V(h)(Y)=\int_{Q_T}(-\Delta)^{\sigma/2} Z(Y,Y')h(Y')\,dY',
\\[10pt]
\label{eq:Vt}
&\partial_t V(h)(Y)=h(Y)+\int_{Q_T}a(Y')(-\Delta)^{\sigma/2} Z(Y,Y')h(Y')\,dY',
\end{align}
and, as a consequence,
\begin{equation}
\label{eq:LV}
L V(h)(Y)=h(Y)+\int_{Q_T}LZ(Y,Y')h(Y')\,dY'.
\end{equation}
\end{Lemma}
\noindent{\it Proof. } In order to prove \eqref{eq:Lambda V} we only have to show that the integral on the right-hand side is well defined. The convergence of the integral at infinity is clear, as always, by the decay of the functions and the fact that $h$ belongs to $L_\sigma$. To see the convergence at the origin we use the representation $(-\Delta)^{\sigma/2} Z(Y,Y')=(-\Delta)^{\sigma/2} W(Y,Y',a(Y'))$ defined before. But we observe that this last function satisfies, for every $Y\in Q_T$, $0<\varepsilon<R$, and $s\neq0$,
$$
\int_{\varepsilon<|A(s)(Y-Y')|_\sigma<R}(-\Delta)^{\sigma/2} W(Y,Y',s)\,dY'=\frac1{|s|}\int_{\varepsilon<|\eta|_\sigma<R}(-\Delta)^{\sigma/2} P(\eta)\,d\eta=0,
$$
by Proposition~2.2 in \cite{VPQR}. Thus, putting $\Omega=\{\varepsilon<|A(a(Y))(Y-Y')|<R\}$, we have
$$
\begin{array}{rl}
\displaystyle\left|\int_\Omega\right.&\displaystyle\left.(-\Delta)^{\sigma/2} Z(Y,Y')h(Y')\,dY'\right|\\ [4mm]
&\displaystyle\le\int_\Omega|(-\Delta)^{\sigma/2} W(Y,Y',a(Y'))-(-\Delta)^{\sigma/2} W(Y,Y',a(Y))|\,|h(Y')|\,dY' \\ [4mm]
&\quad+\displaystyle\int_\Omega|(-\Delta)^{\sigma/2} W(Y,Y',a(Y))|\,|h(Y)-h(Y')|\,dY'.
\end{array}
$$
The first term is estimated using the Mean Value Theorem,
$$
|(-\Delta)^{\sigma/2} W(Y,Y',a(Y'))-(-\Delta)^{\sigma/2} W(Y,Y',a(Y))|\le \frac{c}{|Y-Y'|_\sigma^{N+\sigma-\alpha}}.
$$
The second term is estimated easily by the regularity of $h$ and \eqref{eq:est-lambdaZ}.
Define now $J(x,t,\tau)=\int_{\mathbb{R}^N}Z(x,t,\xi,\tau)f(\xi,\tau)\,d\xi$. It satisfies
$$
L_0J=0\;\text{ for every } 0<\tau<t,\qquad \lim_{\tau\to t}J(x,t,\tau)=f(x,t).
$$
Then, computing
$$
\partial_tV(x,t)=J(x,t,t)+\int_0^T\partial_tJ(x,t,\tau)\,d\tau
$$
we get~\eqref{eq:Vt}.
 \qed

In order to finish the proof of Theorem~\ref{th:parametrix} we still have to show that $L\Gamma=0$ and that $\Gamma$ takes a Dirac delta as initial value, giving a justification to the formal calculus \eqref{formalLGamma}. First we divide the integral in \eqref{Gamma} in two time intervals, $[\tau,t_1]\cup[t_1,t]$. Then we apply Lemma~\ref{lem-volumepot}, with $h(Y)$ equal to the function $\Phi(Y,\oy)$ in \eqref{Gamma(Y)}, to the integral in $[t_1,t]$. Since the function $\Phi$ is H\"older continuous, \eqref{eq:LV} holds. On the other hand, in the interval $[\tau,t_1]$ we use that $(-\Delta)^{\sigma/2}Z$ is continuous and $\Phi$ is absolutely integrable.
Finally, the fact that $\lim_{t\to\tau}\Gamma(x,t,\xi,\tau)=\delta(x-\xi)$ follows immediately by checking that the last integral in \eqref{Gamma(Y)} is convergent, which is now easy using all the estimates obtained so far.
\qed

\noindent{\it Proof of Theorem \ref{th:cauchy}. }
The function $f$ defined through the representation formula~\eqref{sol-linear} solves the problem, all the terms appearing in the equation are continuous, and takes the initial datum in a continuous way. Uniqueness can be proved by the well-known method of the adjoint problem
\begin{equation}
\label{eq:linear-adjoint}
\left\{
\begin{array}{ll}
\partial_t\Gamma^*-(-\Delta)^{\sigma/2} (a(x,t)\Gamma^*)=0, & \qquad  x\in\mathbb{R}^N,\; T<t<\tau,
\\ [4mm]
\Gamma^*(x,\tau,\xi,\tau) = \delta(x-\xi), & \qquad x\in\mathbb{R}^N;
\end{array}
\right.
\end{equation}
see again \cite{Friedman}. It is easy to see that the function $\Gamma^*(Y,\oy)=\Gamma(\oy,Y)$ is its fundamental solution.

The $L^\infty$--norm is estimated by the maximum principle. The other two properties are obtained by multiplying the equation by $(-\Delta)^{\sigma/2}f$ and integrating in space and time.
\qed

\vskip 1cm


\noindent{\large \textbf{Acknowledgments}}

\noindent All authors supported by the Spanish project MTM2014-53037-P.

\vskip 1cm




\

\noindent\textbf{Addresses:}

\noindent\textsc{A. de Pablo: } Departamento de Matem\'{a}ticas, Universidad
Carlos III de Madrid, 28911 Legan\'{e}s, Spain. (e-mail: arturo.depablo@uc3m.es).

\noindent\textsc{F. Quir\'{o}s: } Departamento de Matem\'{a}ticas, Universidad
Aut\'{o}noma de Madrid, 28049 Madrid, Spain. (e-mail: fernando.quiros@uam.es).

\noindent\textsc{A. Rodr\'{\i}guez: } Departamento de Matem\'{a}tica, ETS
Arquitectura, Universidad Polit\'{e}cnica de Madrid, 28040 Madrid, Spain.
(e-mail: ana.rodriguez@upm.es).


\begin{thebibliography}{99}                                                                                               %

\bibitem{Andreu-Mazon-Rossi-Toledo} Andreu, F.; Maz\'on, J.\,M.; Rossi, J.\,D.; Toledo, J. \emph{A nonlocal p-Laplacian evolution equation with nonhomogeneous Dirichlet boundary conditions.} SIAM J. Math. Anal. 40 (2008/09), no.\,5, 1815--1851.

\bibitem{Aronson-Serrin} Aronson, D.\,G.; Serrin, J. \emph{Local behavior of
solutions of quasilinear parabolic equations.} Arch. Rational Mech. Anal. 25
(1967), 81--122.

\bibitem{Athanasopoulos-Caffarelli} Athanasopoulos, I.; Caffarelli, L.\,A.
\emph{Continuity of the temperature in boundary heat control problems.}  Adv.
Math. 224 (2010), no.\,1, 293--315.

\bibitem{Barlow-Bass-Chen-Kassmann} Barlow, M.\,T.; Bass, R.\,F.; Chen, Z.-Q.; Kassmann, M. \emph{Non-local Dirichlet forms and symmetric jump processes.} Trans. Amer. Math. Soc. 361 (2009), no.\,4, 1963--1999.

\bibitem{Bass-Kassmann-Kumagai} Bass, R.\,F.; Kassmann, M.; Kumagai, T. \emph{Symmetric jump processes: localization, heat kernels and convergence.} Ann. Inst. Henri Poincaré Probab. Stat. 46 (2010), no.\,1, 59--71.

\bibitem{Biler-Dolbeault-Esteban} Biler, P.; Dolbeault, J.; Esteban, M.\,J. \emph{Intermediate asymptotics in L1 for general nonlinear diffusion equations.} Appl. Math. Lett. 15 (2002), no.\,1, 101--107.



\bibitem{Biler-Imbert-Karch-2015} Biler, P.; Imbert, C.; Karch, G. \emph{The nonlocal porous medium equation: Barenblatt profiles and other weak solutions.} Arch. Ration. Mech. Anal. 215 (2015), no.\,2, 497--529.

\bibitem{Blumenthal-Getoor} Blumenthal, R.\,M.; Getoor, R.\,K. \emph{Some
theorems on stable processes.} Trans. Amer. Math. Soc. 95 (1960), no.\,2, 263--273.

\bibitem{Bonforte-Figalli-RosOton} Bonforte, M.; Figalli, A.; Ros-Oton, X. \emph{Infinite speed of propagation and regularity of solutions
to the fractional porous medium equation
in general domains.} Comm. Pure Appl. Math., to appear.  \texttt{arXiv:1510.03758v1 [math.AP]}. 


\bibitem{Brandle-dePablo} Br\"andle, C.; de Pablo, A.; \emph{Nonlocal heat equations:
decay estimates and Nash inequalities.} Preprint, \texttt{arXiv:1312.4661v3 [math.AP]}.


\bibitem{Caffarelli-Chan-Vasseur} Caffarelli, L.; Chan, C.\,H.; Vasseur, A. \emph{Regularity
theory for parabolic nonlinear integral operators.} J. Amer. Math. Soc. 24
(2011), no.\,3, 849--869.

\bibitem{Caffarelli-Friedman} Caffarelli, L.\,A.; Friedman, A. \emph{Regularity of the free boundary of a gas flow in an n-dimensional porous medium.} Indiana Univ. Math. J. 29 (1980), no.\,3, 361--391.

\bibitem{Caffarelli-Silvestre} Caffarelli, L.; Silvestre, L. \emph{An
extension problem related to the fractional Laplacian.} {Comm. Partial
Differential Equations} 32 (2007), no.\,7-9, 1245--1260.

\bibitem{Caffarelli-Soria-Vazquez} Caffarelli, L.; Soria, F.; V\'azquez, J.\,L. \emph{Regularity of solutions of the fractional porous medium flow.} J. Eur. Math. Soc. 15  (2013),  no.\,5, 1701--1746.

\bibitem{Caffarelli-Vazquez-ARMA-2011} Caffarelli, L.; V\'azquez, J.\,L. \emph{Nonlinear porous medium flow with fractional potential pressure.} Arch. Ration. Mech. Anal. 202 (2011), no.\,2, 537--565.

\bibitem{Caffarelli-Vazquez-DCDS-2011} Caffarelli, L.\,A.; V\'azquez, J.\,L. \emph{Asymptotic behaviour of a porous medium equation with fractional diffusion.} Discrete Contin. Dyn. Syst. 29 (2011), no.\,4, 1393--1404.


\bibitem{Carrillo-diFrancesco-Toscani} Carrillo, J.\,A.; Di Francesco, M.; Toscani, G. \emph{Intermediate asymptotics beyond homogeneity and self-similarity: long time behavior for $u_t =\Delta\phi(u)$.} Arch. Ration. Mech. Anal.  180  (2006),  no.\,1, 127--149.

\bibitem{ChangLara-Davila} Chang-Lara, H.; D\'avila, G. \emph{Regularity for solutions of non local parabolic equations.} Calc. Var. 49 (2014), 139--172.

\bibitem{Chen} Chen, Z.-Q. \emph{Symmetric jump processes and their heat kernel estimates.} Sci. China Ser. A 52 (2009), no.\,7, 1423--1445.

\bibitem{Chen-Zhang} Chen, Z.-Q.; Zhang, X. \emph{Heat kernels and analyticity of non-symmetric jump diffusion semigroups.} Probab. Theory Relat. Fields, to appear.   DOI:  10.1007/s00440-015-0631-y.

\bibitem{Crandall-Liggett} Crandall, M.\,G.; Liggett, T.\,M. \emph{Generation of
semi-groups of nonlinear transformations on general Banach spaces}. Amer. J.
Math. 93 (1971), 265--298.

\bibitem{deGiorgi} De Giorgi, E. \emph{Sulla differenziabilit\'a e
l'analiticit\'a delle estremali degli integrali multipli regolari}.
(Italian) Mem. Accad. Sci. Torino. Cl. Sci. Fis. Mat. Nat. (3) 3 (1957) 25--43.

\bibitem{PQRV} de Pablo, A.; Quir\'os, F.; Rodr\'{\i}guez, A.; V\'azquez,
J.\,L.  \emph{A fractional porous medium equation.} Adv. Math. 226 (2011),
no.\,2, 1378--1409.

\bibitem{PQRV2} de Pablo, A.; Quir\'os, F.; Rodr\'{\i}guez, A.; V\'azquez,
J.\,L.  \emph{A general fractional porous medium equation.}  Comm. Pure Appl.
Math. 65 (2012), no.\,9, 1242--1284.

\bibitem{PQRV3} de Pablo, A.; Quir\'os, F.; Rodr\'{\i}guez, A.; V\'azquez,
J.\,L.  \emph{Classical solutions for a logarithmic fractional diffusion
equation.}  J. Math. Pures Appl. (9) 101 (2014), no.\,6, 901--924.

\bibitem{dePablo-Vazquez} de Pablo, A.; V\'azquez, J.\,L. \emph{Regularity of solutions and interfaces of a generalized porous medium equation in $R^N$.} Ann. Mat. Pura Appl. (4) 158 (1991), 51--74.

\bibitem{Endal-Jakobsen-delTeso} Endal, J.; Jakobsen, E.\,R.;  del Teso, F. \emph{Uniqueness and properties of distributional solutions of nonlocal degenerate diffusion equations of porous medium type.} Preprint, \texttt{arXiv:1507.04659 [math.AP]}.

\bibitem{Felsinger-Kassmann} Felsinger, M.; Kassmann, M. \emph{Local regularity for parabolic nonlocal operators.} Comm. Partial
Differential Equations 38 (2013), no.\,9, 1539--1573.

\bibitem{Friedman} Friedman, A. \lq\lq Partial differential equations of
parabolic type''. Prentice-Hall, Inc., Englewood Cliffs, N.\,J., 1964.

\bibitem{Hardy-Littlewood} Hardy, G.\,H.; Littlewood, J.\, E. \emph{Some
properties of fractional integrals.~I}. Math. Z. 27 (1928), no.\,1, 565--606.

\bibitem{Kamin}
Kamin, S. \emph{Similar solutions and the asymptotics of filtration equations.}
Arch. Rational Mech. Anal.  60  (1975/76), no.\,2, 171--183.

\bibitem{Kamin-Vazquez} Kamin, S.; V\'azquez, J.\,L. \emph{Asymptotic behaviour of solutions of the porous medium equation with changing sign.} SIAM J. Math. Anal. 22 (1991), no.\,1, 34--45.

\bibitem{Kassmann-2009} Kassmann, M. \emph{A priori estimates for integro-differential operators with measurable kernels.} Calc. Var. Partial Differential Equations 34 (2009), no.\,1, 1--21.

\bibitem{Komatsu}  Komatsu, T. \emph{Uniform estimates for fundamental solutions associated with non-local Dirichlet forms.} Osaka J. Math. 32 (1995), no.\,4, 833--860.

\bibitem{Levi} Levi, E.\,E. \emph{Sulle equazioni lineari totalmente
ellittiche alle derivate parziali.}  Rend. Circ. Mat. Palermo 24 (1907), no.\,1, 275--317 .

\bibitem{Moser} Moser, J. \emph{A Harnack inequality for parabolic differential equations.} Comm. Pure Appl. Math. 17 (1964), 101--134.

\bibitem{Moser-1971} Moser, J. \emph{On a pointwise estimate for parabolic differential equations.} Comm. Pure Appl. Math. 24 (1971), 727--740.

\bibitem{Oleinik-Kalashnikov-Czou} Oleinik, O.\,A.; Kalashnikov, A.\,S.;  Czou, Y.-I.
    \emph{The Cauchy problem and boundary problems for equations of the type of non-stationary filtration.}
    Izv. Akad. Nauk SSSR. Ser. Mat. 22 (1958), 667--704 (Russian).

\bibitem{Schilling-Uemura} Schilling, R.\,L.; Uemura, T. \emph{On the Feller property of Dirichlet forms generated by pseudo differential operators.} Tohoku Math. J. (2) 59 (2007), no.\,3, 401--422.

\bibitem{Serra} Serra, J. \emph{Regularity for fully nonlinear nonlocal parabolic equations with rough kernels.}   Calc. Var. Partial Differential Equations 54 (2015), no.\,1, 615--629.

\bibitem{Sobolev} Sobolev, S.\,L. \emph{On a theorem of functional analysis}.
Transl. Amer. Math. Soc. 34(2) (1963), 39--68; translation of Mat. Sb. 4 (1938) 471--497.

\bibitem{Stan-delTeso-Vazquez-2015} Stan, D.; del Teso, F.; V\'azquez, J.\,L.\emph{Transformations of self-similar solutions for porous medium equations of fractional type.} Nonlinear Anal. 119 (2015), 62--73.

\bibitem{Varopoulos} Varopoulos, N.\,T. \emph{Hardy-{L}ittlewood theory for
semigroups}. J. Funct.  Anal. 63 (1985), no.\,2, 240--260.

\bibitem{Vazquez} V\'{a}zquez, J.\,L. \emph{Barenblatt solutions and
asymptotic behaviour for a nonlinear fractional heat equation of porous medium
type.} J. Eur. Math. Soc. 16 (2014), no.\,4, 769--803.

\bibitem{VPQR} V\'azquez, J.\,L.; de Pablo, A.; Quir\'os, F.; Rodr\'{\i}guez,
A. \emph{Classical solutions and higher regularity for nonlinear fractional
diffusion equations.} J. Eur. Math. Soc., to appear. \texttt{arXiv:1311.7427 [math.AP]}.
\end{thebibliography}
\end{document}